\documentclass[10pt,a4paper]{amsart}
\usepackage{amssymb}
\newtheorem{theorem}{Theorem}[section]
\newtheorem{proposition}[theorem]{Proposition}
\newtheorem{lemma}[theorem]{Lemma}
\newtheorem{corollary}[theorem]{Corollary}

\theoremstyle{remark}
\newtheorem*{remark}{Remark}
\newtheorem*{remarks}{Remarks}
\newtheorem*{example}{Example}

\newcommand{\set}[1]{\left\{#1\right\}}
\def\uu{\mathcal U}
\def\vv{\mathcal V}
\def\uuu{\overline{\mathcal U}}
\def\NN{\mathbb N}
\numberwithin{equation}{section}

\begin{document}

\title{Unique expansions of real numbers}
\author{Martijn de Vries}
\address{Delft University of Technology, Mekelweg 4, 2628 CD Delft, the Netherlands}
\email{w.m.devries@tudelft.nl}
\author{Vilmos Komornik}

\address{D\'epartement de Math\'ematique,
         Universit\'e Louis Pasteur,
         7 rue Ren\'e Descartes, 67084 Strasbourg Cedex, France}
\email{komornik@math.u-strasbg.fr}
\thanks{The first author has been supported by NWO Project nr. ISK04G}
\dedicatory{Dedicated to the memory of Paul Erd\H os}
\begin{abstract}
It was discovered some years ago that there exist non-integer real numbers
$q>1$ for which only one sequence $(c_i)$ of integers $c_i \in [0,q)$ satisfies
the equality $\sum_{i=1}^\infty c_iq^{-i}=1$. The set of such
``univoque numbers'' has a rich topological structure, and its study revealed
a number of unexpected connections with measure theory, fractals, ergodic theory
and Diophantine approximation.

In this paper we consider for each fixed $q>1$ the set
$\uu_q$ of real numbers $x$
having a unique representation of the form $\sum_{i=1}^\infty
c_iq^{-i}=x$ with integers $c_i$ belonging to $[0,q)$. We carry out a detailed
topological study of these sets. For instance, we characterize
their closures, and we determine those bases $q$ for which $\uu_q$
is closed or even a Cantor set. We also study the set $\uu_q'$ consisting of
all sequences $(c_i)$ of integers $c_i \in [0,q)$ such that $\sum_{i=1}^{\infty} c_i q^{-i} \in \uu_q$. We determine the numbers $r >1$ for which the map $q \mapsto \uu_q'$ (defined on $(1, \infty)$) is constant in a neighborhood of $r$ and the numbers $q >1$ for which $\uu_q'$
is a subshift or a subshift of finite type.
\end{abstract}

\maketitle

\section{Introduction and statement of the main results}\label{s1}

Following a seminal paper of R\'enyi \cite{[R]} many works were devoted
to probabilistic, measure theoretical and number theoretical
aspects of representations for real numbers in non-integer bases; see, e.g., Frougny
and Solomyak \cite{[FS]}, Peth\H o and Tichy \cite{[PT]}, Schmidt \cite{[Sc]}. A new
research field was opened when
 Erd\H os, Horv\'ath and Jo\'o \cite{[EHJ]} discovered continuum many real numbers
$q>1$ for which only one sequence $(c_i)=c_1 c_2 \ldots$ of integers $c_i$ belonging to $[0,q)$ satisfies
the equality
\begin{equation*}
\sum_{i=1}^{\infty} \frac{c_i}{q^i}=1.
\end{equation*}
(They considered the case $1<q<2$.)
Subsequently, the set ${\uu}$ of all such {\em univoque numbers} $q>1$ was
characterized lexicographically in \cite{[EJK1], [KL3]}, its smallest element was determined in [KL1], and its
topological structure was described in \cite{[KL3]}. On the other hand, the
investigation of numbers $q>1$ for which there exist continuum many
such sequences,
including sequences containing all possible finite variations of the
integers $c \in [0,q)$, revealed close connections to Diophantine approximations;
see, e.g.,
\cite{[EJK1], [EJK3], [EK], [KLP]},  Borwein and Hare \cite{[BH1], [BH2]}, Komatsu \cite{[K]},
 and Sidorov \cite{[Si1]}.

For any fixed real number $q>1$, we may also introduce the set ${\uu}_q$
of real numbers $x$ for which exactly one sequence $(c_i)$ of integers $c_i \in [0,q)$ satisfies the equality
\begin{equation*}
\sum_{i=1}^{\infty} \frac{c_i}{q^i}=x.
\end{equation*}
If $q$ is an integer, these sets are well-known.
However, their structure is more complex if $q$ is a non-integer,
see, e.g., Dar\'oczy and K\'atai \cite{[DK1], [DK2]}, Glendinning and Sidorov \cite{[GS]},
and Kall\'os \cite{[K1], [K2]}. The purpose of this paper is to give a complete
topological description of the sets ${\uu}_q$: they have a different nature for
different classes of the numbers  $q$. Our investigations also provide new
results concerning the set $\uu$ of univoque numbers. For instance, we determine for each $n \in \NN:=\set{1,2, \ldots}$ the smallest element of $\uu \cap (n,n+1)$ and we continue the study of the topological structure of $\uu$, started in \cite{[KL3]}.
In order to state our results we need to introduce some notation and terminology.

In this paper a {\it sequence} always means an element of the set $\set{0,1,2, \ldots}^{\NN}$. A sequence is called {\it infinite} if it contains infinitely many nonzero elements; otherwise it is called {\it finite}. We use systematically the lexicographical order between sequences: we write $(a_i) < (b_i)$ or $(b_i) > (a_i)$ if there exists an index $n \in \NN$ such that $a_i=b_i$ for $i < n$ and $a_n < b_n$. We also equip for each $n \in \NN$ the set $\set{0,1,2, \ldots}^{n}$ of \textit{blocks of length} $n$ with the lexicographical order.

Given a real number $q > 1$, an {\it expansion in base} $q$ (or simply {\it expansion}) of
a real number $x$ is a sequence $(c_i)$ such that
\begin{equation*}
0 \leq c_i < q  \mbox{ for all } i \geq 1, \quad \mbox{and} \quad
x=\sum_{i=1}^{\infty} \frac{c_i}{q^i}.
\end{equation*}
If a real number $x$ has an expansion in base $q$, then $x$ must belong to the interval
\begin{equation*}
J_q := \Bigl[0,\frac{\lceil q \rceil -1}{q-1}\Bigr]
\end{equation*}
where $\lceil q \rceil$ denotes the smallest integer larger than or equal to $q$. Note that
$1 \in J_q$.

A \textit{sequence} $(c_i)$ such that $0 \leq c_i < q$ for all $i \geq 1$ is called
\textit{univoque} in base $q$ if
\begin{equation*}
x=\sum_{i=1}^{\infty}\frac{c_i}{q^i}
\end{equation*}
is an element of ${\uu}_q$.

The \textit{greedy} expansion
$(b_i(x,q))=(b_i(x))=(b_i)$ of a number
$x \in J_q$ in base $q$
is the largest expansion of $x$ in lexicographical order. It is
well-known that the greedy expansion of any $x \in J_q$
exists; \cite{[R],[P],[EJK2]}.
A \textit{sequence} $(b_i)$ is called \textit{greedy} in base $q$
if $(b_i)$ is the greedy expansion of
\begin{equation*}
x=\sum_{i=1}^{\infty} \frac{b_i}{q^i}.
\end{equation*}

The \textit{quasi-greedy} expansion $(a_i(x,q))=(a_i(x))=(a_i)$ of a
number $x \in J_q \setminus \set{0}$ in base $q$ is the largest
{\it infinite} expansion
of $x$ in lexicographical order. Observe that we have to
exclude the number $0$ since there do not exist infinite expansions
of $x=0$ at all. On the other hand, the largest infinite expansion
of any $x \in J_q \setminus \set{0}$ exists, as we shall prove in the
next section. In order to simplify some statements below, the quasi-greedy
expansion of the number $0 \in J_q$ is defined to be $0^{\infty}=00 \ldots$.
Note that this is the only expansion of $x=0$.
A \textit{sequence} $(a_i)$ is called \textit{quasi-greedy} in base $q$
if $(a_i)$ is the quasi-greedy expansion of
\begin{equation*}
x=\sum_{i=1}^{\infty} \frac{a_i}{q^i}.
\end{equation*}
We usually denote the quasi-greedy expansion $(a_i(1,q))$ of the number
$1$ in base $q$ by $(\alpha_i(q))=(\alpha_i)$. To stress that the quasi-greedy expansion of 1 in base $q$ is given by $(\alpha_i)$, we sometimes
write $q \sim (\alpha_i)$. This notation is particularly convenient in Section~\ref{s6} where we consider expansions $(\alpha_i(q))$ for different values of $q$ simultaneously.

Since $\alpha_1=\lceil q \rceil -1$ (as we shall see in the next section), the digits $c_i$ of an expansion
$(c_i)$ belong to $A:=\set{0,\ldots , \alpha_1}$ for all $i \ge 1$.
Hence we consider expansions with coefficients or digits in the
alphabet $A$ of numbers $x \in [0, \alpha_1/(q-1)]$.

The greedy expansion of a number
$x \in J_q \setminus \set{0}$ coincides with the quasi-greedy expansion
if and only if the greedy expansion of $x$ is infinite.
If the greedy expansion $(b_i)$ of $x \in J_q \setminus \set{0}$ is finite and
$b_n$ is its last
nonzero element, then the quasi-greedy expansion of $x$ is given by
\begin{equation*}
(a_i)=b_1 \ldots b_{n-1} b_n^- \alpha_1 \alpha_2 \ldots, \quad \text{where} \quad b_n^-:=b_n-1.
\end{equation*}
For instance, if $q$ equals the {\it golden ratio} $G:=(1+\sqrt{5})/2$ and
$x=q^{-1} + q^{-2} + q^{-3}$, then $(\alpha_i)=(10)^{\infty}$,
$(b_i(x))= 1110^{\infty}$, and $(a_i(x))=1 (10)^{\infty}$.

Of course, whether a sequence is univoque, greedy or quasi-greedy
depends on the base $q$.
However, when $q$ is understood from the context, we simply speak of univoque sequences and
(quasi)-greedy sequences.
Furthermore, we shall write
$\overline{c}:=\alpha_1 -c \, (c \in A)$, unless stated otherwise.
We shall also write $\overline{c_1 \ldots c_n}$ instead of
$\overline{c_1} \ldots \overline{c_n}$ and
$\overline{c_1 c_2 \ldots}$ instead of $\overline{c_1} \,
\overline{c_2} \ldots$ $(c_i \in A, i \ge 1)$.
Sometimes we refer to $\overline{c_1 c_2 \ldots}$ as the {\it conjugate}
of an expansion $(c_i)$. Finally, we set $c^+:=c+1$ ($c \in A \setminus \set{\alpha_1}$) and $c^-:=c-1$ ($c \in A \setminus \set{0}$).

The following important theorem, which is essentially due to Parry \cite{[P]}(see also \cite{[DK1], [DK2]}),
 plays a crucial role in the proof of our main results:
\begin{theorem}\label{t11}
Fix $q>1$.
\begin{itemize}
\item[\rm (i)]
A sequence $(b_i)=b_1b_2 \ldots \in \set{0, \ldots, \alpha_1}^{\NN}$
is greedy if and only if
\begin{equation*}
{b_{n+1}b_{n+2}\ldots} <
\alpha_1 \alpha_2\ldots \quad \mbox{whenever} \quad  b_n < \alpha_1.
\end{equation*}
\item[\rm (ii)] A sequence $(c_i)=c_1c_2\ldots \in
\set{0,\ldots,\alpha_1}^{\NN}$
is univoque if and only if
\begin{equation*}
c_{n+1}c_{n+2}\ldots < \alpha_1 \alpha_2 \ldots \quad \mbox{whenever} \quad
c_n < \alpha_1
\end{equation*}
and
\begin{equation*}
\overline{c_{n+1}c_{n+2}\ldots} < \alpha_1 \alpha_2 \ldots \quad
\mbox{whenever} \quad  c_n >0.
\end{equation*}
\end{itemize}
\end{theorem}
\noindent
Note that $(\alpha_i)$ is the unique expansion of $1$ in base $q$
if and only if $q \in \uu$. Hence, replacing the sequence $(c_i)$
in Theorem~\ref{t11}~(ii) with the sequence $(\alpha_i)$, one
obtains a characterization of $\uu$.

Recently, the authors of \cite{[KL3]} studied the topological
structure of the set ${\uu}$. In particular, they
showed that ${\uu}$ is {\it not} closed and they characterized
its closure $\uuu$:

\begin{theorem} \label{t12}
A real number $q>1$ belongs to $\uuu$ if and only if the quasi-greedy
expansion $(\alpha_i)$ of the number $1$ in base $q$ satisfies
\begin{equation*}
\overline{\alpha_{k+1}\alpha_{k+2}\ldots }< \alpha_1\alpha_2\ldots  \quad
\mbox
{for all} \quad k \geq 1.
\end{equation*}
\end{theorem}
\noindent
It is possible to give a similar description of the set $\uuu$ in words: $q > 1$  belongs to $\uuu$ if and only if $(\alpha_i(q))$
is the unique {\it infinite} expansion of the number 1 in base $q$ (see Corollary~\ref{c54}).

\begin{remarks} \mbox{}
\begin{itemize}
\item Recall that the number $q$ is not allowed in any expansion in base $q$ if $q$ is an integer.
Our choice of the digit set simplifies some statements.
For example it will follow from the theorems below that
\begin{equation*}
\uu_q = \overline{\uu_q} \quad \Longleftrightarrow \quad q \in (1, \infty)
\setminus \uuu
\end{equation*}
where $\overline{\uu_q}$ denotes the closure of $\uu_q$.
\item In the definition of $\uu$ given in \cite{[KL3]} the integers $2,3,\ldots$ were excluded.
However, its closure is the same as the set $\uuu$ defined in this paper. As a consequence, Theorem~\ref{t12} still holds in our set-up.
\item It follows from Theorems~\ref{t11},
\ref{t12} and Proposition~\ref{p23} below  that the quasi-greedy expansion of 1 in base $q$ is periodic for each $q \in \uuu \setminus \uu$. For instance, if $n \ge 2$, then $(1^n0)^{\infty}= (\alpha_i(q))$ for some $ q \in \uuu \setminus \uu$.
\end{itemize}
\end{remarks}

For any fixed $q>1$, we introduce the set $\vv_q$, defined by
\begin{equation*}
{\vv}_q= \set{x \in J_q : \overline{a_{n+1}(x) a_{n+2}(x) \ldots}
\leq \alpha_1 \alpha_2 \ldots \quad \mbox{whenever} \quad a_n(x) >
0 }.
\end{equation*}
It follows from Theorem~\ref{t11} that $\uu_q \subset \vv_q$ for
all $q >1$. \

Now we are ready to state our main results.

\begin{theorem} \label{t13}
Suppose that $q \in {\uuu}$. Then
\begin{itemize}

\item[\rm (i)] $\overline{\uu_q} = \vv_q$;
\item[\rm (ii)] $\vert \vv_q \setminus {\uu}_q \vert = \aleph_0$
and $\vv_q \setminus {\uu}_q$ is dense in $\vv_q$;
\item[\rm (iii)] if $q \in {\uu}$, then each element $x \in \vv_q
\setminus \uu_q$ has exactly $2$ expansions;
\item[\rm (iv)] if $q \in {\uuu} \setminus {\uu}$, then each
element $x \in \vv_q \setminus \uu_q$ has exactly $\aleph_0$ expansions.
\end{itemize}
\end{theorem}

\begin{remarks}\mbox{}

\begin{itemize}
\item The proof of part~(ii) yields the following more precise
results where for $q \in {\uuu}$ we set
\begin{equation*}
A_q=\set{ x \in \vv_q \setminus \uu_q : x \mbox{ has a finite
greedy expansion}}
\end{equation*}
and
\begin{equation*}
B_q=\set{ x \in \vv_q \setminus \uu_q : x \mbox{ has an infinite
greedy expansion}}:
\end{equation*}
\begin{itemize}
\item If $q \in {\uuu}\setminus \NN$, then both $A_q$ and $B_q$
are countably infinite and dense in $\vv_q$. Moreover, the greedy
expansion of a number $x \in B_q$ ends with
$\overline{\alpha_1\alpha_2\ldots}$, and
\begin{equation*}
B_q = \set{ \alpha_1/(q-1) - x : x \in A_q}.
\end{equation*}
\item If $q \in \set{2,3,\ldots}$, then $B_q=\varnothing$.
\end{itemize}
\item For each $x \in \vv_q \setminus {\uu}_q$, the proof of parts
(iii) and (iv) also provides the list of all expansions of $x$ in
terms of its greedy expansion.
\end{itemize}
\end{remarks}

Our next goal is to describe the relationship between the sets ${\uu}_q$,
$\overline{{\uu}_q}$ and $\vv_q$ in case $q \notin {\uuu}$. To
this end, we introduce the set ${\vv}$, consisting of those
numbers $q >1$, for which the quasi-greedy expansion $(\alpha_i)$ of the number
1 in base $q$ satisfies
\begin{equation*}
\overline{\alpha_{k+1} \alpha_{k+2} \ldots} \leq \alpha_1 \alpha_2
\ldots \quad \mbox{for all} \quad k \geq 1.
\end{equation*}
It follows from Theorem~\ref{t12} that ${\uu} \subset {\uuu}
\subset {\vv}$. The following results combined with Theorem~\ref{t13} imply that ${\uu}_q$ is
closed if $q \notin {\uuu}$ and that the set ${\vv}_q$ is closed
for all $q >1$.

\begin{theorem}\label{t14}
Suppose that $q \in {\vv} \setminus {\uuu}$. Then
\begin{itemize}
\item[\rm (i)]
the sets ${\uu}_q$ and ${\vv}_q$ are closed;
\item[\rm (ii)]
$\vert {\vv}_q \setminus {\uu}_q \vert = \aleph_0$ and
${\vv}_q \setminus {\uu}_q$ is a discrete set, dense in
${\vv}_q$;
\item[\rm (iii)]
each element $x \in {\vv}_q \setminus {\uu}_q$ has exactly $\aleph_0$
expansions and a finite greedy expansion.
\end{itemize}
\end{theorem}

\begin{remark}
Our proof also provides the list of all
expansions of all elements $x \in {\vv}_q \setminus
{\uu}_q$.
\end{remark}

\begin{theorem}\label{t15}
Suppose that $q \in (1, \infty) \setminus {\vv}$. Then
\begin{equation*}
{\uu}_q = \overline{{\uu}_q}= {\vv}_q.
\end{equation*}
\end{theorem}

\begin{remarks}\mbox{}

\begin{itemize}
\item In view of the above results,
Theorem~\ref{t11}
already gives us a lexicographical characterization of $\overline{{\uu}_q}$ if
$q \in (1, \infty) \setminus \uuu$ because in this case ${\uu}_q$ is closed.
\item It is well-known that the set ${\uu}$ has Lebesgue measure
zero; \cite{[EHJ],[KK]}. In \cite{[KL3]} it was shown that the set ${\uuu}
\setminus {\uu}$ is countably infinite. It follows from the above
results that ${\uu}_q$ is closed for almost every $q >1$.
\item Let $q>1$ be a non-integer.
In \cite{[DDV]} it has been proved that almost
every $x \in J_q$ has a continuum of expansions in base $q$ (see also [Si2]).
It follows from the above results that the set
$\overline{{\uu}_q}$ has Lebesgue measure zero. Hence the set
${\uu}_q$ is nowhere dense.
\item Let  $q >1$ be an integer.
In this case the quasi-greedy expansion of
$1$ in base $q$ is given by $(\alpha_i)=\alpha_1^{\infty}=(q-1)^{\infty}.$
Moreover, the set
$J_q \setminus {\uu}_q$ is countably infinite
and each element in $J_q \setminus {\uu}_q$ has only two expansions,
one of them being finite while the other one ends with an infinite string of
$(q-1)$'s.
\item In \cite{[KL1]} it was shown that the smallest element of $\uu$ is
given by $q'\approx 1.787$, and the unique expansion of 1 in base
$q'$ is the truncated Thue--Morse sequence
$(\tau_i)=\tau_1 \tau_2 \ldots = 11010011 \ldots$, which can be defined recursively
by setting $\tau_{2^N}=1$ for $N=0,1,2,\ldots$ and
\begin{equation*}
\tau_{2^N+i} = 1 - \tau_i \quad \mbox{for }
1 \leq i < 2^N, \, N=1,2, \ldots.
\end{equation*}
Subsequently, Glendinning and Sidorov \cite{[GS]} proved that ${\uu}_q$
is countable \footnote{Here and in the sequel, we call a set \textit{countable} if it is either finite or countably infinite.}
if $1< q < q'$ and has the cardinality of the
continuum if $q \in [q',2)$. Moreover, they showed that ${\uu}_q$
is a set of positive Hausdorff dimension if $q' < q < 2$, and they
described a method to compute its Hausdorff dimension (see also
\cite{[DK2], [K1], [K2]}).
\end{itemize}
\end{remarks}

In the following theorem we characterize those $q>1$ for which ${\uu}_q$ or
$\overline{{\uu}_q}$ is a Cantor set, i.e., a nonempty closed set having
neither interior nor isolated points. We recall from \cite{[KL3]} that

\begin{itemize}
\item $\vv$ is closed and $\uu$ is closed from above \footnote{We call a set $X \subset \mathbb{R}$ \textit{closed from above}
(\textit{closed from below})
if for each $x \in \mathbb{R} \setminus X$ there exists a number $\delta >0$ such that
$[x,x+\delta) \cap X = \varnothing$ ($(x-\delta, x] \cap X = \varnothing$).},
\item $\vert\uuu\setminus\uu\vert=\aleph_0$ and $\uuu\setminus\uu$ is dense
in $\uuu$,
\item $\vert\vv\setminus\uuu\vert=\aleph_0$ and $\vv\setminus\uuu$ is a
discrete set, dense in $\vv$.
\end{itemize}

Since the set $(1, \infty) \setminus \vv$ is open, we can write
$(1,\infty) \setminus \vv$ as the union of countably many disjoint open
intervals $(q_1,q_2)$: its connected components. Let us denote by $L$ and $R$ the set of
left (respectively right) endpoints of the intervals $(q_1,q_2)$.

\begin{theorem}\label{t16}\mbox{}
\begin{itemize}
\item[\rm (i)] $L={\NN}\cup(\vv\setminus\uu)$ and $R=\vv\setminus\uuu$.
Hence $R \subset L$, and
\begin{equation*}
 (1,\infty)\setminus\uuu=\cup (q_1,q_2]
\end{equation*}
where the union runs over the connected components $(q_1,q_2)$ of
$(1,\infty)\setminus\vv$.
\item[\rm (ii)] If $q \in \set{2,3, \ldots}$, then neither $\uu_q$ nor
$\overline{{\uu}_q}$ is a Cantor set.
\item[\rm (iii)] If $q\in\uuu \setminus \NN$, then ${\uu}_q$ is not a
Cantor set, but its closure $\overline{{\uu}_q}$ is a Cantor set.
\item[\rm (iv)] If $q\in (q_1,q_2]$, where $(q_1,q_2)$ is a
connected component of $(1,\infty)\setminus\vv$, then the closed set
${\uu}_q$ is
a Cantor set if and only if $q_1\in \set{3,4, \ldots} \cup
({\uuu}\setminus \uu)$. Moreover, if $q_1 \in \set{1,2} \cup (\vv \setminus \uuu)$, then the isolated
points of $\uu_q$ form a dense subset of $\uu_q$.
\end{itemize}
\end{theorem}
\begin{remark}
We also describe the set of endpoints of the connected
components $(p_1,p_2)$ of $(1, \infty) \setminus \uuu$: denoting by $L'$
and $R'$ the set of left (respectively right) endpoints of the intervals $(p_1,p_2)$,
we have
\begin{equation*}
L' = {\NN} \cup (\uuu \setminus \uu)
\quad\mbox{and}\quad   R' \subset \uu.
\end{equation*}
This enables us to determine the condensation points of $\uu_q$ for each $q >1$;
see the remarks at the end of Section~\ref{s6}.
\end{remark}

The ideas leading to the above theorem result in a new characterization of
{\em stable bases}, introduced and investigated by
Dar\'oczy and K\'atai (\cite{[DK1],[DK2]}).  Let us denote by $\uu_q'$ and
$\vv_q'$ the sets of quasi-greedy expansions in base $q$ of all
numbers $x\in\uu_q$ and $x\in\vv_q$ respectively. Note that $\uu_q'$ is simply the set of univoque
sequences in base $q$.

If $1 < q < r$ and $\uu_q' = \uu_r'$, then $\lceil q \rceil = \lceil r \rceil$ and
$\uu_q' = \uu_t'$ for each $t \in (q,r)$, as follows from Theorem~\ref{t11}
and Proposition~\ref{p23} below. For this reason, we call a number $q >1$
{\em stable from above} (respectively {\em stable from below}) if there
exists a number $s>q$ (respectively $1< s<q$) such that
\begin{equation*}
\uu_q'= \uu_s'.
\end{equation*}
We call a number $q>1$ \textit{stable} if it is stable from above and from below.
Finally, we say that an interval $I\subset (1,\infty)$ is a
{\em stability interval} if $\uu_q'=\uu_s'$ for all $q,s\in I$.

\begin{theorem} \label{t17}
The maximal stability intervals are given by the singletons $\set{q}$ where
$q \in \uuu$ and the intervals $(q_1,q_2]$ where
$(q_1,q_2)$ is a connected component of $(1,\infty)\setminus\vv$. Moreover, if
$q_1 \in \vv \setminus \uu$, then
\begin{equation*}
\uu_q'= \vv_{q_1}'\quad\text{for all}\quad q\in (q_1,q_2].
\end{equation*}
\end{theorem}
\begin{remark}
The proof of Theorem~\ref{t17} yields a new characterization of
the sets $\uuu$ and $\vv$ (see Proposition~\ref{p69}).
\end{remark}

We recall (see, e.g., \cite{[LM]}) that a set $S \subset \set{0, \ldots,
\alpha_1}^{\NN}$ is called a {\it subshift} if
there exists a set $\mathcal{F}(S) \subset \cup_{k=1}^{\infty}
\set{0, \ldots, \alpha_1}^{k}$ such that a sequence $(c_i) \in
\set{0, \ldots, \alpha_1}^{\NN}$ belongs to $S$ if and only if none of the blocks $c_{i+1} \ldots c_{i+n}$ $(i \ge 0, n \ge 1)$ belongs to $\mathcal{F}(S)$. A subshift $S$ is called a
{\it subshift of finite type} if one can choose $\mathcal{F}(S)$ to be finite. We endow the set $\set{0, \ldots, \alpha_1}^{\NN}$
with the topology of coordinate-wise convergence.

\begin{theorem}\label{t18}
Let $q >1$ be a real number. The following statements are equivalent.
\begin{itemize}
\item[\rm(i)] $q \in (1, \infty) \setminus \uuu$.
\item[\rm(ii)] $\uu_q'$ is a subshift of finite type.
\item[\rm(iii)] $\uu_q'$ is a subshift.
\item[\rm(iv)] $\uu_q'$ is a closed subset of $\set{0, \ldots, \alpha_1}^{\NN}$.
\end{itemize}
\end{theorem}

Finally, we determine the cardinality of $\uu_q$ for all $q>1$.
We recall that for $q \in (1,2)$
this has already been done by Glendinning and Sidorov (\cite{[GS]}),
using a different method. Let $q''$ be the smallest element of
$\uu \cap (2,3)$.
It was shown in \cite{[KL2]} that the unique expansion of $1$ in base $q''$ is given
by $(\lambda_i)=\lambda_1 \lambda_2 \ldots=21020121 \ldots$, where
$\lambda_i = \tau_i + \tau_{2i -1}$, $ i \ge 1$.

\begin{theorem}\label{t19}
Let $q >1$ be a real number.\mbox{}
\begin{itemize}
\item[\rm (i)] If $q \in (1, G]$, then $\uu_q$ consists merely of the endpoints of $J_q$.
\item[\rm (ii)] If $q \in (G, q') \cup (2, q'')$, then
$\vert \uu_q \vert = \aleph_0$.
\item[\rm (iii)] If $q \in [q',2] \cup [q'', \infty)$, then $\vert \uu_q
\vert = 2^{\aleph_0}.$
\end{itemize}
\end{theorem}

\begin{remarks}\mbox{}
\begin{itemize}
\item We also determine the unique expansion of $1$ in base $q^{(n)}$
for $n \in \set{3,4,\ldots}$ where $q^{(n)}$ denotes the smallest
element of $\uu \cap (n, n+1)$; see the remarks at the end of
Section~\ref{s6}.
\item A generalization of Theorem~\ref{t19} can be found in \cite{[DV2]}.
\end{itemize}
\end{remarks}

For the reader's convenience we recall some properties of
quasi-greedy expansions in the next section. These properties are
also stated in \cite{[BK]} and are closely related to some important
results, first established in the seminal works by R\'enyi \cite{[R]} and
Parry \cite{[P]}. In Section~\ref{s3} we derive some preliminary lemmas
needed for the proof of our main results. Section~\ref{s4} is
then devoted to  the proof of Theorem~\ref{t13}.
Theorems~\ref{t14} and \ref{t15} are proved in Section~\ref{s5},
and our final Theorems~\ref{t16}, ~\ref{t17}, ~\ref{t18} and
~\ref{t19} are established in Section~\ref{s6}.

\section{Quasi-greedy expansions}\label{s2}

Let $q>1$ be a real number and let $m = \lceil q \rceil -1$.
In the previous section we defined the quasi-greedy expansion of $x \in (0, m /(q-1)]$ as its
largest infinite expansion in base $q$.
In order to prove that this notion is well-defined, we introduce the
{\it quasi-greedy algorithm}: if for some $n \in \NN$,
$a_i(x)=a_i$ is already defined for $i$ with $1 \le i <n$ (no condition if $n=1$), then $a_n(x)=a_n$ is the largest
element of the digit set $A=\set{0, \ldots , m }$ such that
\begin{equation*}
\sum_{i=1}^{n}\frac{a_i}{q^i} < x.
\end{equation*}
Of course, this definition is only meaningful if $x >0$.
In the following proposition we show that this algorithm generates an
expansion of $x$ for all $x \in J_q \setminus \set{0}$. It follows that the
quasi-greedy expansion of $x \in J_q \setminus \set{0}$ is obtained by performing the quasi-greedy algorithm.

\begin{proposition}\label{p21}
Let $x \in (0, m /(q-1)]$. Then
\begin{equation*}
x=\sum_{i=1}^{\infty}\frac{a_i}{q^i}.
\end{equation*}
\end{proposition}

\begin{proof}
If $x = m /(q-1)$, then the quasi-greedy algorithm
provides
$a_i= m $ for all $i \geq 1$ and the desired equality follows.

Suppose that $x \in (0, m /(q-1))$. Then, by definition of the
quasi-greedy
algorithm, there exists an index $n$ such that $a_n < m$.

First assume that $a_n < m$ for infinitely many $n$. For any such $n$,
we have by definition
\begin{equation*}
0 < x-\sum_{i=1}^{n}\frac{a_i}{q^i} \leq \frac{1}{q^n}.
\end{equation*}
Letting $n \to \infty$, we obtain
\begin{equation*}
x=\sum_{i=1}^{\infty}\frac{a_i}{q^i}.
\end{equation*}

Next assume there exists a largest $n$ such that $a_n < m $. Then
\begin{equation*}
\sum_{i=1}^{n}\frac{a_i}{q^i} + \sum_{i=n+1}^{N}\frac{m}{q^i} < x
\leq \sum_{i=1}^{n}\frac{a_i}{q^i}+\frac{1}{q^n},
\end{equation*}
for each $N>n$. Hence
\begin{equation*}
\sum_{i=n+1}^{\infty}\frac{m}{q^i}
\leq x-\sum_{i=1}^{n}\frac{a_i}{q^i}
\leq \frac{1}{q^n}.
\end{equation*}
Note that
\begin{equation*}
\frac{1}{q^n} \leq \sum_{i=n+1}^{\infty}\frac{m}{q^i}
\end{equation*}
for any $q>1$, and
\begin{equation*}
\frac{1}{q^n} = \sum_{i=n+1}^{\infty}\frac{m}{q^i}
\end{equation*}
if and only if $q= m+1$. Hence $q$ is an integer, and
\begin{equation*}
x=\sum_{i=1}^{n}\frac{a_i}{q^i} + \sum_{i=n+1}^{\infty} \frac{m}
{q^i}= \sum_{i=1}^{\infty}\frac{a_i}{q^i}.\qedhere
\end{equation*}
\end{proof}

Now we consider the quasi-greedy expansion $(\alpha_i(q))=(\alpha_i)$ of $x=1$.
Note that $\alpha_1 = m = \lceil q \rceil - 1$ by definition of the quasi-greedy algorithm.

\begin{lemma}\label{l22}
For each $n \geq1$, the inequality
\begin{equation}\label{21}
\sum_{i=1}^{n} \frac{\alpha_i}{q^i} + \frac{1}{q^n} \geq 1
\end{equation}
holds.
\end{lemma}

\begin{proof}
The proof is by induction on $n$.
For $n=1$, the inequality holds because $\alpha_1 +1 \geq q$.
Assume the inequality is valid for some $n \in {\NN}$.
If $\alpha_{n+1} < \alpha_1$, then \eqref{21} with $n+1$ instead of $n$ follows from the
definition of the quasi-greedy algorithm. If $\alpha_{n+1}=\alpha_1$, then
the same conclusion follows from the induction hypothesis and the inequality
$\alpha_1 +1 \geq q$.
\end{proof}

\begin{proposition}\label{p23}
The map $q \mapsto (\alpha_i(q))$ is a strictly increasing bijection from
the open interval $(1, \infty)$ onto the set of all {\rm infinite} sequences $(\alpha_i)$,
satisfying
\begin{equation}\label{22}
\alpha_{k+1} \alpha_{k+2} \ldots \leq \alpha_1 \alpha_2 \ldots \quad
\mbox{for all} \quad k\geq 1.
\end{equation}
\end{proposition}

\begin{proof}
By definition of the quasi-greedy algorithm, the map $q \mapsto(\alpha_i(q))$ is
strictly increasing. Fix $q >1$ and $k \in \NN$. By the preceding lemma we have
\begin{equation*}
\sum_{i=1}^{k} \frac{\alpha_i(q)}{q^i} +\frac{1}{q^k}
\geq \sum_{i=1}^{\infty} \frac{\alpha_i(q)}{q^i}
\end{equation*}
whence
\begin{equation}\label{23}
\sum_{i=1}^{\infty} \frac{\alpha_{k+i}(q)}{q^i} \le 1.
\end{equation}
If we had $(\alpha_{k+i}(q)) > (\alpha_i(q))$, then by Lemma~\ref{l22},
\begin{equation*}
\sum_{i=1}^{n} \frac{\alpha_{k+i}(q)}{q^i} \ge 1
\end{equation*}
for some $n \in \NN$, which contradicts \eqref{23} because $(\alpha_{k+i}(q))$ is infinite.

Conversely, let $(\alpha_i)$ be an infinite sequence
satisfying \eqref{22}. Solving the equation
\begin{equation*}
\sum_{i=1}^{\infty}\frac{\alpha_i}{q^i}=1,
\end{equation*}
we obtain a unique number $q >1$. Note that $0 \le \alpha_n \le \alpha_1 < q$ for $n \ge 1$. In order to prove that $(\alpha_i)=(\alpha_i(q))$,
it suffices to show that for each
$n\ge 1$, the inequality
\begin{equation*}
\sum_{i=n+1}^{\infty} \frac{\alpha_i}{q^i} \leq \frac{1}{q^n}
\end{equation*}
holds. Starting with $k_0:=n$ and using \eqref{22}, we {\it try} to define a
sequence
\begin{equation*}
k_0<k_1<\cdots
\end{equation*}
satisfying for $j=1,2,\ldots$ the conditions
\begin{equation*}
\alpha_{k_{j-1}+i}=\alpha_i \quad \mbox{for }i=1,\ldots,k_j-k_{j-1}-1
\quad \mbox{and } \alpha_{k_j} < \alpha_{k_j-k_{j-1}}.
\end{equation*}
If we obtain in this way an infinite number of indices, then we have
\begin{align*}
\sum_{i=n+1}^{\infty} \frac{\alpha_i}{q^i}&\leq
\sum_{j=1}^{\infty} \left( \left(\sum_{i=1}^{k_j-k_{j-1}}
\frac{\alpha_i}{q^{k_{j-1}+i}} \right)-\frac{1}{q^{k_j}} \right)\\
&<
\sum_{j=1}^{\infty} \left(\frac{1}{q^{k_{j-1}}}-\frac{1}{q^{k_j}} \right)=
\frac{1}{q^n}.
\end{align*}
If we only obtain a finite number of indices, then there exists a least
nonnegative integer $N$ ($N=0$ is possible) such that
$(\alpha_{k_N+i})=(\alpha_i)$ and we have
\begin{align*}
\sum_{i=n+1}^{\infty} \frac{\alpha_i}{q^i}&\leq
\sum_{j=1}^{N} \left(\left(\sum_{i=1}^{k_j-k_{j-1}}
\frac{\alpha_i}{q^{k_{j-1}+i}}\right)-\frac{1}{q^{k_j}} \right)
+\sum_{i=1}^{\infty} \frac{\alpha_i}{q^{k_N+i}}\\
&\leq
\sum_{j=1}^{N} \left(\frac{1}{q^{k_{j-1}}}-\frac{1}{q^{k_j}} \right)
+\sum_{i=1}^{\infty} \frac{\alpha_i}{q^{k_N+i}} =
\frac{1}{q^n}. \qedhere
\end{align*}
\end{proof}

The proof of the following propositions is almost identical to the proof
of Proposition~\ref{p23} and is therefore omitted.

\begin{proposition}\label{p24}
For each $q>1$, the map $x \mapsto(a_i(x))$ is a strictly increasing bijection
from $(0, \alpha_1 /(q-1)]$ onto the set of all {\rm infinite} sequences $(a_i)$,
satisfying
\begin{equation*}
0 \leq a_n \leq \alpha_1 \quad \mbox{for all} \quad n \geq 1
\end{equation*}
and
\begin{equation}\label{24}
a_{n+1}a_{n+2}\ldots \leq \alpha_1 \alpha_2 \ldots \quad \mbox{whenever}
\quad a_n < \alpha_1.
\end{equation}
\end{proposition}

\begin{proposition}\label{p25}
For each $q>1$, the map $x \mapsto (b_i(x))$ is a strictly increasing bijection
from $[0, \alpha_1/(q-1)]$ onto the set of all sequences $(b_i)$, satisfying
\begin{equation*}
0 \leq b_n \leq \alpha_1 \quad \mbox{for all} \quad n \geq 1
\end{equation*}
and
\begin{equation}\label{25}
b_{n+1}b_{n+2} \ldots < \alpha_1 \alpha_2 \ldots \quad \mbox{whenever}
\quad b_n < \alpha_1.
\end{equation}
\end{proposition}

\begin{remarks}\mbox{}
\begin{itemize}
\item A sequence $(c_i)$ is univoque if and
only if $(c_i)$ is greedy and $(\alpha_1 - c_i)$ is greedy. Hence
Theorem~\ref{t11} is a consequence of Proposition~\ref{p25}.
\item The greedy expansion $(b_i)$ of $x \in J_q$ is generated
by the {\it greedy algorithm}: if for some $n \in \NN$, $b_i$
is already defined for $i$ with $1 \le i < n$ (no condition if $n=1$),
then $b_n$ is the largest element of $A$ such that
\begin{equation*}
\sum_{i=1}^{n} \frac{b_i}{q^i} \leq x.
\end{equation*}
The proof of this assertion goes along the same lines as the proof of
Proposition~\ref{p21}.
\end{itemize}
\end{remarks}

\section{Some preliminary results}\label{s3}

Throughout this section, $q >1 $ is an arbitrary but fixed real number.

\begin{lemma}\label{l31}
Let $(d_i)=d_1d_2 \ldots$ be a greedy or quasi-greedy sequence.
Then for all $N \geq 1$ the truncated
sequence $d_1 \ldots d_N 0^{\infty}$ is greedy.
\end{lemma}

\begin{proof}
If $(d_i)=0^{\infty}$, then there is nothing to prove. If $(d_i)\not=0^{\infty}$, then the statement follows from
Propositions~\ref{p24} and \ref{p25}.
\end{proof}

\begin{lemma}\label{l32}
Let $(b_i) \not=\alpha_1^{\infty}$ be a greedy sequence and let $N$ be a
positive integer. Then there exists a greedy sequence $(c_i)> (b_i)$ such
that
\begin{equation*}
c_1 \ldots c_N= b_1\ldots b_N.
\end{equation*}
\end{lemma}

\begin{proof}
Since $(b_i) \not=\alpha_1^{\infty}$, it follows from
\eqref{25} that $b_n < \alpha_1$
for infinitely many $n$. Hence we may assume, by enlarging $N$ if
necessary, that $b_N<\alpha_1$. Let
\begin{equation*}
I=\set{i \in \NN: 1\le i \le N \text{ and }  b_i < \alpha_1}=:\set{i_1,\ldots,i_n}.
\end{equation*}
Note that for $i_r \in I$,
\begin{equation*}
\sum_{j=1}^{\infty} \frac{b_{i_r+j}}{q^j}=
\sum_{j=1}^{N-i_r}\frac{b_{i_r+j}}{q^j} + \frac{1}{q^{N-i_r}}
\sum_{i=1}^{\infty} \frac{b_{N+i}}{q^i}<1
\end{equation*}
because $(b_i)$ is greedy and $b_{i_r}< \alpha_1$. For each $r \in \set{1, \ldots, n}$, choose $y_r$ such that
\begin{equation}\label{31}
\sum_{i=1}^{\infty} \frac{b_{N+i}}{q^i} < y_r \leq \frac{\alpha_1}{q-1}
\end{equation}
and
\begin{equation}\label{32}
\sum_{j=1}^{N-i_r}\frac{b_{i_r+j}}{q^j} + \frac{y_r}{q^{N-i_r}}
<1.
\end{equation}
Let $y=\min\set{y_1,\ldots,y_n}$ and denote the greedy expansion of
$y$ by $d_1d_2\ldots$. Finally, let $(c_i)=b_1\ldots b_N d_1d_2\ldots$. From
\eqref{31} we infer that $(c_i)>(b_i)$. It remains to show
that $(c_i)$ is a greedy sequence, i.e., we need to show that
\begin{equation}\label{33}
\sum_{i=1}^{\infty} \frac{c_{j+i}}{q^i} < 1 \quad \mbox{whenever} \quad
c_j < \alpha_1.
\end{equation}
If $c_j < \alpha_1$ and $j \leq N$, then \eqref{33}
follows from \eqref{32}. If $c_j < \alpha_1$ and $j>N$, then
\eqref{33} follows from the
fact that $(d_i)$ is a greedy sequence.
\end{proof}

\begin{lemma}\label{l33}
If $(b_i) \not= \alpha_1^{\infty}$ is a greedy sequence, then there
exists a sequence $1 \leq n_1 < n_2 < \cdots$ such that for each
$i \geq 1$,
\begin{equation}\label{34}
b_{n_i} < \alpha_1 \quad \mbox{and} \quad
b_{m+1} \ldots b_{n_i} < \alpha_1 \ldots \alpha_{n_i-m} \quad
\mbox{if }  m <n_i  \text{ and }  b_m < \alpha_1.
\end{equation}
\end{lemma}

\begin{proof}
We define a sequence $(n_i)_{i \geq 1}$ satisfying the
requirements by induction.

Let $r$ be the least positive integer for which $b_r < \alpha_1$.
Then, \eqref{34} with $r$ instead of $n_i$ holds clearly. Set $n_1:=r$ and let $\ell$ be a positive
integer.

Suppose we have already defined
$n_1 < \cdots < n_{\ell}$ such that \eqref{34} holds for each $i$ with $1 \le i \le \ell$.
Since $(b_i)$ is greedy and $b_{n_{\ell}}< \alpha_1$,
there exists a {\it smallest} integer
$n_{\ell+1} > n_{\ell}$ such that
\begin{equation}\label{35}
b_{n_{\ell}+1} \ldots b_{n_{\ell+1}} <  \alpha_1 \ldots \alpha_{n_{\ell+1}
- n_{\ell}}.
\end{equation}
Note that $b_{n_{\ell + 1}} < \alpha_{n_{\ell+1} - n_{\ell}} \leq \alpha_1$.
It remains to verify that
\begin{equation}\label{36}
b_{m+1} \ldots b_{n_{\ell +1}} < \alpha _1 \ldots \alpha_{n_{\ell + 1} - m}
\end{equation}
if $ 1 \le m < n_{\ell+1}$ and $b_m < \alpha_1$. If $m < n_{\ell}$, then \eqref{36}
follows from the induction hypothesis. If $m = n_{\ell}$, then \eqref{36} reduces
to \eqref{35}. If $n_{\ell} <m < n_{\ell+1}$, then
\begin{equation*}
b_{n_{\ell} + 1} \ldots b_m = \alpha_1 \ldots \alpha_{m- n_{\ell}},
\end{equation*}
by minimality of $n_{\ell +1}$, and thus by \eqref{35} and \eqref{22},
\begin{equation*}
b_{m+1} \ldots b_{n_{\ell+1}} <  \alpha_{m-n_{\ell} + 1} \ldots
\alpha_{n_{\ell+1} - n_{\ell}}
\le  \alpha_1 \ldots \alpha_{n_{\ell + 1} - m}. \qedhere
\end{equation*}
\end{proof}

We call a set $B \subset J_q$ {\it symmetric} if $\ell(B)=B$, where $\ell: J_q \to J_q$ is given by
\begin{equation*}
\ell(x)=\alpha_1 / (q-1) - x \quad \quad  [x \in J_q].
\end{equation*}

\begin{lemma}\label{l34}
\mbox{}
\begin{itemize}
\item[\rm (i)] The sets $\uu_q$ and $\vv_q$ are symmetric.
\item[\rm (ii)] The set $\vv_q$ is closed.
\end{itemize}
\end{lemma}

\begin{proof}
(i) The set $\uu_q$ is symmetric because $(c_i)$ is an expansion of $x$ if and only if
$(\alpha_1 - c_i)$ is an expansion of $\ell(x)$.

If $q$ is a non-integer and $x \in J_q$, then by \eqref{24}, the sequence $(\alpha_1 - a_i(x))$ is either infinite or
is equal to $0^{\infty}$. It follows from Proposition~\ref{p24} that the set $\vv_q$ is symmetric and
$(a_i(\ell(x)))= (\alpha_1 - a_i(x))$ for each $x\in\vv_q$. If $q$ is an integer, then $\vv_q = J_q = [0,1]$.

(ii) We prove that $\vv_q$ is closed by showing that its complement is
open. If $(a_i)$ is the quasi-greedy expansion of some $x \in J_q
\setminus \vv_q$, then there exists an integer $n >0$ such that
\begin{equation*}
a_n >0 \quad \mbox{and} \quad \overline{a_{n+1} a_{n+2} \ldots}
> \alpha_1 \alpha_2 \ldots .
\end{equation*}
Let $m$ be such that
\begin{equation}\label{37}
\overline{a_{n+1} \ldots a_{n+m}} > \alpha_1 \ldots \alpha_m,
\end{equation}
and let
\begin{equation*}
y=\sum_{i=1}^{n+m} \frac{a_i}{q^i}.
\end{equation*}
According to Lemma~\ref{l31} the greedy expansion of $y$ is given by
$a_{1} \ldots a_{n+m}0^{\infty}$. Therefore the quasi-greedy expansion of each
number $v \in (y,x]$ starts with the block $a_1 \ldots a_{n+m}$.
It follows from \eqref{37} that
\begin{equation*}
(y,x] \cap {\vv}_q = \varnothing.
\end{equation*}
Since $x \in J_q \setminus \vv_q$ is arbitrary and $\vv_q$ is symmetric, there also exists a number $z > x$ such that
\begin{equation*}
[x, z) \cap \vv_q = \varnothing.
\qedhere
\end{equation*}
\end{proof}

\begin{lemma}\label{l35}
Let $(b_i)$ be the greedy expansion of some $x \in [0, \alpha_1/(q-1)]$ and
suppose that for some $n \geq 1$,
\begin{equation*}
b_n > 0 \quad \mbox{and} \quad \overline{b_{n+1}b_{n+2}\ldots} > \alpha_1\alpha_2 \ldots .
\end{equation*}
Then
\begin{itemize}
\item[\rm (i)] there exists a number $z > x$ such that $[x,z] \cap
{\uu}_q =\varnothing$;
\item[\rm (ii)] if $b_j >0$ for some $j >n$, then there exists a number
$y < x$ such
that $[y,x] \cap {\uu}_q = \varnothing$.
\end{itemize}
\end{lemma}

\begin{proof}
(i) Choose a positive integer $M >n$ such that
\begin{equation*}
\overline{b_{n+1}\ldots b_M} > \alpha_1 \ldots \alpha_{M-n}.
\end{equation*}
Applying Lemma~\ref{l32} choose a greedy sequence $(c_i) > (b_i)$
such that $c_1\ldots c_M$ $=b_1 \ldots b_M$. Then $(c_i)$ is the greedy
expansion of some $z > x$.
If $(d_i)$ is the greedy expansion of some element in $[x,z]$, then
$(d_i)$ also begins with  $b_1\ldots b_M$ and hence
\begin{equation*}
\overline{d_{n+1}\ldots d_M}> \alpha_1 \ldots \alpha_{M-n}.
\end{equation*}
In particular, we have
\begin{equation*}
d_n>0
\quad\text{and}\quad
\overline{d_{n+1}d_{n+2} \ldots}> \alpha_1\alpha_2 \ldots .
\end{equation*}
We infer from Theorem~\ref{t11} that
$[x,z] \cap {\uu}_q = \varnothing.$

(ii) Suppose that $b_j >0$ for some $j > n$. It follows from Lemma~\ref{l31} that
$(c_i):=b_1 \ldots b_n 0^{\infty}$ is the greedy expansion of some $y < x$.
If $(d_i)$ is the greedy expansion
of some element in $[y,x]$, then $(c_i) \leq (d_i) \leq (b_i)$ and
$d_1 \ldots d_n=b_1 \ldots b_n$.
Therefore
\begin{equation*}
\overline{d_{n+1}d_{n+2}\ldots} \geq \overline{b_{n+1}b_{n+2}\ldots}
> \alpha_1\alpha_2 \ldots,
\end{equation*}
and $d_n=b_n >0$. It follows from Theorem~\ref{t11} that
$[y,x] \cap {\uu}_q = \varnothing$.
\end{proof}

\section{Proof of Theorem~\ref{t13}}\label{s4}

If $q$ belongs to $\uuu$, then we know from Theorem~\ref{t12} and
Proposition~\ref{p23} that the quasi-greedy expansion $(\alpha_i)$ of 1 in base
$q$ satisfies
\begin{equation}\label{41}
\alpha_{k+1} \alpha_{k+2} \ldots \leq \alpha_1 \alpha_2 \ldots \quad
\mbox{for all} \quad  k \geq 1
\end{equation}
and
\begin{equation}\label{42}
\overline{\alpha_{k+1} \alpha_{k+2} \ldots } < \alpha_1 \alpha_2 \ldots \quad
\mbox{for all} \quad  k \geq 1.
\end{equation}
Note that a sequence $(\alpha_i)$ satisfying \eqref{41} and \eqref{42} is
automatically infinite, and is thus the quasi-greedy expansion of 1 in base $q$ for some
$q \in \uuu$.
The following lemmas are obtained in \cite{[KL3]}.

\begin{lemma}\label{l41}
If $(\alpha_i)$ is a sequence satisfying \eqref{41} and
\eqref{42}, then there exist arbitrarily large integers $m$ such
that for all $k$ with $0 \le k < m$,
\begin{equation}\label{43}
\overline{\alpha_{k+1} \ldots \alpha_{m}} < \alpha_1 \ldots \alpha_{m-k}.
\end{equation}
\end{lemma}

\begin{lemma}\label{l42}
Let $(\gamma_i)$ be a sequence satisfying
\begin{equation*}
\gamma_{k+1}\gamma_{k+2}\ldots \leq \gamma_1 \gamma_2 \ldots
\end{equation*}
and
\begin{equation*}
\overline{\gamma_{k+1}\gamma_{k+2}\ldots} \leq \gamma_1 \gamma_2 \ldots
\end{equation*}
for all $k \geq 1$, with $\overline{\gamma_j}:=\gamma_1-\gamma_j$, $j \ge 1$.
If
\begin{equation*}
\overline{\gamma_{n+1}\ldots \gamma_{2n}} = \gamma_1 \ldots \gamma_n
\end{equation*}
for some $n \geq 1$, then
\begin{equation*}
(\gamma_i)= (\gamma_1 \ldots \gamma_n \overline{\gamma_1 \ldots \gamma_n})^{\infty}.
\end{equation*}
\end{lemma}

\begin{lemma}\label{l43}
If $q \in \uuu \setminus \uu$, then the greedy
expansion $(\beta_i)$ of $1$ is finite, and all expansions of $1$ are given by
\begin{equation}\label{44}
(\alpha_i) \quad \mbox{and} \quad
 (\alpha_1 \ldots \alpha_m)^N \alpha_1 \ldots \alpha_{m-1} \alpha_m^+
0^{\infty},
\quad N=0,1,2, \ldots,
\end{equation}
where $m$ is such that $\beta_m$ is the last nonzero element of
$(\beta_i)$.
\end{lemma}

\begin{proof}[Proof of Theorem~\ref{t13}]
(i) Fix $q \in \uuu$. It follows from Lemma~\ref{l34} that $\overline{\uu_q} \subset
\vv_q$. Therefore, it suffices to show that each $x \in \vv_q \setminus \uu_q$ belongs to $\overline{\uu_q}$.

First assume that $x \in \vv_q \setminus \uu_q$ has a finite greedy expansion $(b_i)$.
If $b_n$ is the last nonzero element of
$(b_i)$, then
\begin{equation*}
(a_i(x))=(a_i)=b_1\ldots b_n^{-} \alpha_1 \alpha_2 \ldots .
\end{equation*}
According to Lemma~\ref{l41} there exists a sequence
$1 \leq m_1 < m_2 < \cdots $ such that \eqref{43}
is satisfied with $m=m_i$ for all $i \geq 1$.
We may assume that $m_i >n$ for all $i \geq 1$.
Consider for each $i \geq 1$ the sequence $(b_j^i)$, given by
\begin{equation*}
(b_j^i)=b_1 \ldots b_n^-(\alpha_1 \ldots \alpha_{m_i}\overline
{\alpha_1 \ldots \alpha_{m_i}})^{\infty},
\end{equation*}
and define the number  $x_i$ by
\begin{equation*}
x_i= \sum_{j=1}^{\infty}\frac{b_j^i}{q^j}.
\end{equation*}
Note that the sequence $(x_i)_{i \geq 1}$ converges to $x$ as
$i$ goes to infinity. It remains to show that $x_i \in
{\uu}_q$ for all $i \geq 1$.
According to Theorem~\ref{t11} it suffices to verify that
\begin{equation}\label{45}
b_{m+1}^i b_{m+2}^i \ldots < \alpha_1 \alpha_2 \ldots \quad \mbox{whenever}
\quad
b_m^i < \alpha_1
\end{equation}
and
\begin{equation}\label{46}
\overline{b_{m+1}^i b_{m+2}^i \ldots} <
\alpha_1 \alpha_2 \ldots \quad \mbox{whenever} \quad b_m^i > 0.
\end{equation}
According to \eqref{42} we have
\begin{equation*}
\overline{\alpha_{m_i+1} \ldots \alpha_{2m_i}} \leq
\alpha_1 \ldots \alpha_{m_i}.
\end{equation*}
Note that this inequality cannot be an equality, for otherwise it would
follow from Lemma \ref{l42} that
\begin{equation*}
(\alpha_i)=(\alpha_1 \ldots \alpha_{m_i} \overline
{\alpha_1 \ldots \alpha_{m_i}})^{\infty}.
\end{equation*}
However, this sequence does not satisfy \eqref{42} for $k=m_i$.
Therefore
\begin{equation*}
\overline{\alpha_{m_i+1} \ldots \alpha_{2m_i}} <
\alpha_1 \ldots \alpha_{m_i}
\end{equation*}
or equivalently
\begin{equation}\label{47}
\overline{\alpha_1 \ldots \alpha_{m_i}} < \alpha_{m_i+1} \ldots \alpha_{2m_i}.
\end{equation}
If $m \geq n$, then \eqref{45} and \eqref{46}
follow from \eqref{41}, \eqref{43} and \eqref{47}. Now we verify
\eqref{45} and \eqref{46} for $m < n$. Fix $m < n$. If $b_m^i < \alpha_1$, then
\begin{equation*}
b_{m+1}^i \ldots b_n^i = b_{m+1} \ldots b_n^-
< b_{m+1} \ldots b_n
\leq \alpha_1 \ldots \alpha_{n-m},
\end{equation*}
where the last inequality follows from the fact that $(b_i)$ is a
greedy expansion and $b_m=b_m^i < \alpha_1.$ Hence
\begin{equation*}
b_{m+1}^i b_{m+2}^i \ldots < \alpha_1 \alpha_2 \ldots .
\end{equation*}
Suppose that $b_m^i = a_m >0.$ Since
\begin{equation*}
\overline{a_{m+1} a_{m+2} \ldots} \leq \alpha_1 \alpha_2 \ldots
\end{equation*}
by assumption, and $b_{m+1}^i \ldots b_n^i = a_{m+1} \ldots a_n$, it suffices
to verify that
\begin{equation*}
\overline{b_{n+1}^i b_{n+2}^i \ldots} < \alpha_{n-m+1} \alpha_{n-m+2} \ldots .
\end{equation*}
This is equivalent to
\begin{equation} \label{48}
\overline{\alpha_{n-m+1} \alpha_{n-m+2} \ldots} <
 (\alpha_1 \ldots \alpha_{m_i}\overline{\alpha_1 \ldots \alpha_{m_i}})^{\infty}.
\end{equation}
Since $n < m_i$ for all $i \geq 1$, we infer from \eqref{43} that
\begin{equation*}
\overline{\alpha_{n-m+1} \ldots \alpha_{m_i}} < \alpha_1 \ldots
\alpha_{m_i-(n-m)},
\end{equation*}
and \eqref{48} follows.

Next assume that $x \in \vv_q \setminus \uu_q$ has an infinite greedy expansion $(b_i)$.
Since $x \notin {\uu}_q$, there exists a
{\it smallest} positive
integer $n$ such that
\begin{equation}\label{49}
b_n >0 \quad \mbox{and} \quad \overline{b_{n+1}b_{n+2} \ldots}
\geq \alpha_1 \alpha_2 \ldots .
\end{equation}
Since $x \in \vv_q$ and $(a_i(x))=(b_i)$, this last inequality is in fact an
equality.
As before, let $1 \leq m_1 < m_2 < \cdots$ be a sequence such that
\eqref{43} is satisfied with $m=m_i$ for all $i \geq 1$.
Again, we may assume that $m_i >n$ for all $i \geq 1$.
Consider for each $i \geq 1$ the sequence $(b_j^i)$, given by
\begin{equation*}
(b_j^i)=b_1 \ldots b_n (\overline{\alpha_1 \ldots \alpha_{m_i}}
\alpha_1 \ldots \alpha_{m_i})^{\infty},
\end{equation*}
and define the number $x_i$ by
\begin{equation*}
x_i= \sum_{j=1}^{\infty}\frac{b_j^i}{q^j}.
\end{equation*}
Then the sequence $(x_i)_{i \geq 1}$ converges to $x$ as $i$ goes
to infinity. It remains to show that $x_i \in {\uu}_q$ for all $i
\geq 1$, i.e., it remains to verify ~\eqref{45} and ~\eqref{46}.
We leave this easy verification to the reader.

(iia) We establish that $\vert \vv_q \setminus {\uu}_q
\vert = \aleph_0$ for each $q \in \uuu$. More specifically, if $q \in {\uuu}\setminus
{\NN}$, then the sets $A_q$ and $B_q$ (introduced in a remark
following the statement of Theorem~\ref{t13}) are countably
infinite, and the greedy expansion of a number $x \in B_q$
ends with $\overline{\alpha_1\alpha_2\ldots}$. If $q \in
\set{2,3,\ldots }$, then $A_q = \vv_q \setminus \uu_q$.

Fix $q \in {\uuu}$.
Denote the greedy expansion of a number
$x \in \vv_q \setminus {\uu}_q$ by $(b_i)$.
Since $x \notin {\uu}_q$, there exists a number $n$ such that ~\eqref{49} holds.
If both inequalities in ~\eqref{49} are strict, then $b_i = 0$ for all $i >n$,
as follows from Lemma~\ref{l35} and part~(i).
Otherwise, the sequence $(b_i)$ ends with
$\overline{\alpha_1\alpha_2 \ldots}$, which is infinite unless $q$ is an
integer. It follows from Theorems~\ref{t11} and \ref{t12} that a
sequence of the form $0^n10^{\infty}$ ($n \geq 0$) is the
finite greedy expansion of $1/q^{n+1} \in \vv_q
\setminus {\uu}_q$. Moreover, if $q \in {\uuu}
\setminus {\NN}$, then a sequence of the form
$\alpha_1^n \overline{\alpha_1 \alpha_2 \ldots}$
($n \geq 1$) is the infinite
greedy expansion of a number $x \in \vv_q
\setminus {\uu}_q$. These observations conclude the proof.

(iib) We show that if $q \in {\uuu}$, then $A_q$ is dense in
$\vv_q$.

Fix $q \in {\uuu}$. For each $x \in
{\uu}_q$,
we will define a sequence $(x_i)_{i \geq 1}$
of numbers in
$A_q \subset \vv_q \setminus {\uu}_q$ that
converges to $x$. We have seen in the proof of part~(iia) that
$1/q^n \in A_q$ for each $n \geq 1$.
Hence there exists a sequence of numbers in $A_q$
that converges to $0$.
Now suppose that $x \in {\uu}_q \setminus
\set{0}$ and denote by $(c_i)$ the unique expansion of $x$.
Since $\overline{c_1 c_2 \ldots}
\not= \alpha_1^{\infty}$ is
greedy, we infer from Lemma~\ref{l33}
that there exists a sequence $1 \leq n_1 < n_2 < \cdots$, such that
for each $i \geq 1$,
\begin{equation}\label{410}
c_{n_i} > 0 \quad \mbox{and} \quad
\overline{c_{m+1} \ldots c_{n_i}} < \alpha_1
\ldots \alpha_{n_i - m} \quad \mbox{if }  m < n_i \, \,
\mbox{and} \, \, c_m >0.
\end{equation}
Now consider for each $i \geq 1$ the sequence $(b_j^i)$, given by
\begin{equation*}
(b_j^i) = c_1 \ldots c_{n_i} 0 ^{\infty},
\end{equation*}
and define the number $x_i$ by
\begin{equation*}
x_i = \sum_{j=1}^{\infty} \frac{b_j^i}{q^j}.
\end{equation*}
According to Lemma~\ref{l31} the sequence $(b_j^i)$ is the finite
greedy expansion of the number $x_i$, $i \geq 1$. Moreover, the sequence
$(x_i)_{i \geq 1}$ converges to $x$ as $i$ goes to infinity. We claim that
$x_i \in A_q$ for each
$i \geq 1$. Note that $x_i \notin {\uu}_q$ because
the quasi-greedy sequence $(a_j^i)$, given by
\begin{equation*}
c_1 \ldots c_{n_i}^{-} \alpha_1 \alpha_2 \ldots,
\end{equation*}
is another expansion of $x_i$. It remains to prove that
\begin{equation}\label{411}
a_j^i >0    \Longrightarrow \overline{a_{j+1}^i a_{j+2}^i \ldots}
 \leq \alpha_1 \alpha_2 \ldots.
\end{equation}
If $j < n_i$ and $a_j^i>0$, then
\begin{equation*}
\overline{a_{j+1}^i \ldots a_{n_i}^i}
= \overline{c_{j+1} \ldots c_{n_i}^-} \leq \alpha_1 \ldots \alpha_{n_i - j}
\end{equation*}
by \eqref{410}, and
\begin{equation*}
\overline{a_{n_i+1}^i a_{n_i+2}^i \ldots}=
\overline{\alpha_1 \alpha_2 \ldots} <
\alpha_{n_i-j+1} \alpha_{n_i-j+2} \ldots
\end{equation*}
by Theorem~\ref{t12}. If $j=n_i$, then \eqref{411} follows from
$\overline{\alpha_1}=0 < \alpha_1$. Finally, if $j > n_i$, then \eqref{411}
follows again from Theorem~\ref{t12}.

(iic) We show that if $q \in {\uuu} \setminus {\NN}$,
then the set $B_q$ is dense in $\vv_q$, and
\begin{equation*}
B_q = \set{\alpha_1 /(q-1) - x : x \in A_q}.
\end{equation*}

Fix $q \in {\uuu} \setminus {\NN}$ and suppose that $x \in A_q$ has a
finite greedy expansion $(b_i)$ with last nonzero element $b_n$.
An application of Proposition~\ref{p25} and Theorem~\ref{t12} yields that
\begin{equation*}
(c_i)= \overline{b_1 \ldots b_n^- \alpha_1 \alpha_2 \ldots}
\end{equation*}
is the greedy expansion of $\alpha_1 /(q-1) - x$. It follows from the symmetry of $\uu_q$ and $\vv_q$ that
the number $\alpha_1 /(q-1) - x$ belongs to $B_q$.
Conversely, suppose that $x \in B_q$ has an infinite greedy expansion $(b_i)$ and let
$n$ be the {\it smallest} positive integer for which \eqref{49} holds.
Then $b_{n+1}b_{n+2} \ldots = \overline{\alpha_1 \alpha_2 \ldots}$, and
\begin{equation*}
(c_i)= \overline{b_1 \ldots b_n^-} 0 ^{\infty}
\end{equation*}
is the greedy expansion of $\alpha_1 /(q-1) - x \in A_q$.
It follows from the symmetry of $\uu_q$ together with part~(iib)
that the set $B_q$ is dense in $\vv_q$ as well.

(iii) and (iv) Fix $q \in {\uuu}$ and let $(b_i)$ be
the greedy expansion of a number $x \in \vv_q
\setminus \uu_q$. Let $n$ be the {\it smallest} positive integer
for which \eqref{49} holds and let $(d_i)$ be another expansion of $x$. Then
$(d_i) < (b_i)$, and hence there exists a {\it smallest} integer $j\ge1$ for which
$d_j < b_j$. First we show that $j \geq n$. Assume on the contrary that
$j < n$. Then $b_j >0$, and by minimality of $n$ we have
\begin{equation*}
b_{j+1}b_{j+2} \ldots > \overline{\alpha_1 \alpha_2 \ldots}.
\end{equation*}
From Theorems~\ref{t11} and \ref{t12} we know that $\overline{\alpha_1 \alpha_2 \ldots}$ is the greedy expansion of
$\alpha_1/(q-1) - 1$, and thus
\begin{equation*}
\sum_{i=1}^{\infty} \frac{d_{j+i}}{q^i} =
b_j - d_j + \sum_{i=1}^{\infty} \frac{b_{j+i}}{q^i} > 1 + \frac{\alpha_1}{q-1} - 1 =  \frac{\alpha_1}{q-1},
\end{equation*}
which is impossible.
If $j=n$, then $d_n=b_n^-$, for otherwise we have $q > 2$ and
\begin{equation*}
2 \leq \sum_{i=1}^{\infty} \frac{d_{n+i}}{q^i} \leq \frac{\lceil q \rceil - 1}
{q-1},
\end{equation*}
which is also impossible. Now we distinguish between two cases.

If $j=n$ and
\begin{equation}\label{412}
\overline{b_{n+1} b_{n+2} \ldots} > \alpha_1 \alpha_2 \ldots,
\end{equation}
then by Lemma~\ref{l35} and part~(i) we have $b_{r}=0$ for all $r >n$, from which
it follows that $(d_{n+i})$ is an expansion of $1$. Hence, if $q \in
{\uu}$ and \eqref{412} holds, then the only expansion of $x$
starting with $b_1 \ldots b_n^-$ is given by $(c_i):= b_1 \ldots b_n^- \alpha_1
\alpha_2 \ldots$. If $q \in {\uuu} \setminus {\uu}$ and
\eqref{412} holds, then any expansion $(c_i)$ starting with $b_1 \ldots
b_n^-$ is an expansion of $x$ if and only if $(c_{n+i})$ is one of the
expansions listed in
\eqref{44}.

If $j=n$ and
\begin{equation}\label{413}
\overline{b_{n+1} b_{n+2} \ldots} = \alpha_1 \alpha_2 \ldots,
\end{equation}
then
\begin{equation*}
\sum_{i=1}^{\infty} \frac{d_{n+i}}{q^i}=
1+ \sum_{i=1}^{\infty} \frac{b_{n+i}}{q^i}  =
\sum_{i=1}^{\infty} \frac{\alpha_1}{q^i}.
\end{equation*}
Hence, if \eqref{413} holds, then the only expansion of $x$ starting with
$b_1 \ldots b_n^-$ is given by $b_1 \ldots
b_n^- \alpha_1^{\infty}$.

Finally, if $j > n$, then \eqref{413} holds,
for otherwise $(b_{n+i})=0^{\infty}$ and $d_j < b_j$ is impossible.
Note that in this case $ q \notin {\uu}$, because otherwise $(b_{n+i})$ is
the unique expansion of $\sum_{i = 1}^{\infty} \overline{\alpha_i} q^{-i}$
and thus $(d_{n+i}) = (b_{n+i})$ which is impossible due to $j > n$.
Hence, if $q \in {\uu}$, then $(b_i)$ is the only expansion of $x$
starting with $b_1 \ldots b_n$. If $q \in {\uuu}
\setminus {\uu}$ and \eqref{413} holds, then any expansion $(c_i)$
starting with $b_1 \ldots b_n$ is an expansion of $x$ if and only if $(c_{n+i})$
is one of the conjugates of the expansions listed in \eqref{44}.

Parts (iii) and (iv) follow directly from the
above considerations.
\end{proof}

\begin{remarks}
\mbox{}
\begin{itemize}
\item Fix $q \in {\uuu}$. It follows from the proof of Theorem~\ref{t13}~(iia)
that each $x \in \vv_q \setminus {\uu}_q$ has
either a finite expansion or an expansion that ends with
$\overline{\alpha_1 \alpha_2 \ldots}$ in which case $x$ can be written as
\begin{equation*}
x = \frac{b_1}{q} + \cdots + \frac{b_n}{q^n} + \frac{1}{q^n} \left(
\frac{\alpha_1}{q-1} -1 \right).
\end{equation*}
Hence, if $q$ is algebraic, then each $x \in \vv_q
\setminus {\uu}_q$ is algebraic, and if $q$ is transcendental,
then each $x \in \vv_q \setminus {\uu}_q$ is transcendental. If $q \in \uuu \setminus \uu$, then $q$ is algebraic because
$1$ has a finite greedy expansion in base $q$ by Lemma~\ref{l43}.
If $q \in \uu$, then each neighborhood of $q$ contains uncountably many univoque numbers because
$\uuu$ is a perfect set and $\uuu \setminus \uu$ is countable (\cite{[KL3]}). Hence the set of transcendental univoque numbers is dense in $\uuu$. For instance, it was shown by Allouche and Cosnard in \cite{[AC]} that the smallest univoque number $q'$ is transcendental. Subsequently, it was shown in \cite{[DV1]} that the set of algebraic univoque numbers is dense in $\uuu$ as well. This implies in particular that there does not exist a smallest algebraic univoque number, a result first established in \cite{[KLPt]}.
\item It follows from Theorem~\ref{t13}~(ii) and (iii)  that the set
\begin{equation*}
\mathcal{T}_q:= \set{x \in J_q: x \text{ has exactly 2 expansions in base } q}
\end{equation*}
is not closed if $q$ belongs to $\uu$, in which case its closure contains $\uu_q$. It would be interesting to determine all numbers $q >1$ for which $\mathcal{T}_q$ is not closed.
\end{itemize}
\end{remarks}

\section{Proof of Theorems \ref{t14} and \ref{t15}}\label{s5}

Fix $q >1$. It follows from Propositions \ref{p23} and
\ref{p25} that a sequence $(b_i)$ is greedy
if and only if $0 \leq b_n \leq \alpha_1$ for all $ n \geq 1$, and
\begin{equation}\label{51}
b_{n+k+1}b_{n+k+2} \ldots < \alpha_1 \alpha_2 \ldots \quad
\mbox{for all }   k \geq 0, \quad \mbox{whenever} \quad b_n< \alpha_1.
\end{equation}

\begin{lemma}\label{l51}
Suppose that $q \notin {\uuu}$. Then a greedy sequence
$(b_i)$ cannot end with $\overline{\alpha_1\alpha_2 \ldots}$.
\end{lemma}

\begin{proof}
Assume on the contrary that for some $n \ge 0$,
\begin{equation*}
b_{n+1}b_{n+2} \ldots= \overline{\alpha_1 \alpha_2 \ldots}.
\end{equation*}
Since in this case $b_{n+1}=\overline{\alpha_1}=0 < \alpha_1$, it
would follow from \eqref{51} that
\begin{equation*}
\overline{\alpha_{k+1}\alpha_{k+2} \ldots} < \alpha_1 \alpha_2 \ldots
\quad \mbox{for all} \quad k \geq 1.
\end{equation*}
But this contradicts Theorem~\ref{t12}.
\end{proof}

\begin{lemma}\label{l52}
Suppose that $q \notin {\uuu}$. Then

\begin{itemize}
\item[\rm (i)] the set ${\uu}_q$ is closed;
\item[\rm (ii)] each element $x \in {\vv}_q \setminus {\uu}_q$ has
a finite greedy expansion.
\end{itemize}
\end{lemma}

\begin{proof}
(i) Let
$x \in J_q \setminus {\uu}_q$ and denote the
greedy expansion of $x$ in base $q$ by $(b_i)$.
According to Theorem~\ref{t11}
there exists a positive integer $n$ such that
\begin{equation*}
b_n >0 \quad \mbox{and} \quad \overline{b_{n+1}b_{n+2} \ldots} \geq
\alpha_1 \alpha_2 \ldots .
\end{equation*}
Applying Lemmas \ref{l35} and \ref{l51} we conclude
that
\begin{equation*}
[x,z] \cap {\uu}_q = \varnothing
\end{equation*}
for some number $z >x$. It follows that ${\uu}_q$ is closed from above.
Since the set $\uu_q$ is symmetric it is closed from below as well.

(ii) Assume on the contrary that $(a_i(x))=(b_i(x))$ for some $x \in {\vv}_q \setminus {\uu}_q$.
Then it would follow that for some positive integer $n$,
\begin{equation*}
\overline{b_{n+1}(x) b_{n+2}(x) \ldots} = \alpha_1 \alpha_2 \ldots,
\end{equation*}
contradicting Lemma~\ref{l51}.
\end{proof}

Recall from the introduction that the set ${\vv}$ consists
of those numbers $q >1$
for which the quasi-greedy expansion $(\alpha_i)$ of $1$ in base
$q$ satisfies
\begin{equation}\label{52}
\overline{\alpha_{n+1} \alpha_{n+2} \ldots} \leq \alpha_1 \alpha_2 \ldots
\quad \mbox{for all} \quad  n \geq 1.
\end{equation}
If $q \in {\vv} \setminus {\uuu}$, then the quasi-greedy expansion of 1 in base $q$ is of the form
\begin{equation}\label{53}
(\alpha_i)= (\alpha_1 \ldots \alpha_k
\overline{\alpha_1 \ldots \alpha_k})^{\infty},
\end{equation}
where $k$ is the least positive integer satisfying
\begin{equation}\label{54}
\overline{\alpha_{k+1} \alpha_{k+2} \ldots} = \alpha_1 \alpha_2 \ldots.
\end{equation}
In particular, such a sequence is periodic. Note that $\alpha_k >0$, for otherwise it
would follow from \eqref{52} and \eqref{53} that
\begin{equation*}
\overline{\alpha_k \alpha_{k+1}\ldots \alpha_{2k-1}}=\alpha_1 (\alpha_1
\ldots \alpha_{k-1}) \le \alpha_1 \ldots \alpha_{k-1}0,
\end{equation*}
which is impossible because $\alpha_1 >0$ and $\alpha_j \le \alpha_1$ for each $j \in \NN$.
Any sequence of the form $(1^m0^m)^{\infty}$, where $m$ is a positive
integer, is infinite and satisfies \eqref{41} and \eqref{52} but
not \eqref{42}. On the other hand, there are only countably many
periodic sequences. Hence the set ${\vv} \setminus {\uuu}$ is
countably infinite.

The following lemma (\cite{[KL3]}) implies that the number of expansions of 1 is
countably infinite in case $q \in {\vv} \setminus
{\uuu}$. Moreover, all expansions of the number 1 in such a base
$q$ are determined explicitly.

\begin{lemma}\label{l53}
If $q \in {\vv} \setminus {\uuu}$, then all
expansions of $1$ are given by $(\alpha_i)$, and the sequences
\begin{equation*}
(\alpha_1 \ldots \alpha_{2k})^N  \alpha_1 \ldots \alpha_{2k-1} \alpha_{2k}^+
0^{\infty} \, , \quad N=0,1, \ldots
\end{equation*}
and
\begin{equation*}
(\alpha_1 \ldots \alpha_{2k})^N  \alpha_1 \ldots \alpha_{k-1} \alpha_{k}^-
\alpha_1^{\infty} \, , \quad N=0,1, \ldots .
\end{equation*}
\end{lemma}

Now we are ready to prove Theorems \ref{t14} and \ref{t15}.
Throughout the proof of Theorem~\ref{t14}, $q \in {\vv} \setminus
{\uuu}$ is fixed but arbitrary,
 and $k$ is the least positive integer satisfying \eqref{54} with $(\alpha_i)=(\alpha_i(q))$.

\begin{proof}[Proof of Theorem~\ref{t14}]
Thanks to Lemmas \ref{l34} and \ref{l52} we only need to prove parts (ii) and (iii).

(iia) We prove that
$\vert {\vv}_q \setminus {\uu}_q \vert =
\aleph_0$. The set ${\vv}_q \setminus {\uu}_q$ is countable
because
each element $x \in {\vv}_q \setminus {\uu}_q$ has a finite
greedy expansion (see Lemma~\ref{52}). On the other hand, for each
$n \geq 1$ the sequence $\alpha_1^n0^{\infty}$  is the greedy expansion
of an element $x \in {\vv}_q \setminus {\uu}_q$,
from which the claim follows.

(iib) In order to show that ${\vv}_q \setminus {\uu}_q$ is dense
in ${\vv}_q$, one can argue as in the proof of
Theorem~\ref{t13}~(iib). Instead of applying Theorem~\ref{t12} one
should now apply the inequalities \eqref{52}.

(iic) Finally, we show that all elements of ${\vv}_q \setminus {\uu}_q$ are
isolated points of ${\vv}_q$. Let $x \in {\vv}_q \setminus
{\uu}_q$ and let $b_n$ be the last nonzero element of the greedy expansion
$(b_i)$ of $x$.
Choose $m$ such that $\alpha_m < \alpha_1$. This is possible because
$q \notin {\NN}$. According to
Lemma~\ref{l32} there exists a greedy sequence $(c_i) > (b_i)$
such that
\begin{equation*}
c_1 \ldots c_{n+m}=b_1 \ldots b_n 0^{m}.
\end{equation*}
If we set
\begin{equation*}
z= \sum_{i=1}^{\infty} \frac{c_i}{q^i},
\end{equation*}
then the quasi-greedy expansion $(v_i)$ of a number $v \in (x,z]$ starts
with $b_1 \ldots b_n 0^m$.
Hence $v_n = b_n >0$ and
\begin{equation*}
\overline{v_{n+1} \ldots v_{n+m}} = \alpha_1^m > \alpha_1 \ldots \alpha_m.
\end{equation*}
Therefore
\begin{equation*}
(x,z] \cap {\vv}_q = \varnothing.
\end{equation*}
Since the sets $\uu_q$ and $\vv_q$ are symmetric, there also
exists a number $y < x$ satisfying
\begin{equation*}
(y,x) \cap {\vv}_q = \varnothing.
\end{equation*}

(iii) We already know from Lemma~\ref{l52} that each
$x \in {\vv}_q \setminus {\uu}_q$ has a finite greedy expansion.
It remains to show that each element
$x \in {\vv}_q \setminus {\uu}_q$ has exactly $\aleph_0$
expansions.  Let $x \in {\vv}_q \setminus {\uu}_q$ and let
$b_n$ be
the last nonzero element of its greedy expansion $(b_i)$. If $j <n$ and
$b_j = a_j (x) >0$, then
\begin{equation*}
\overline{a_{j+1} (x)\ldots a_n (x)}= \overline{b_{j+1} \ldots b_n^-} \leq
\alpha_1 \ldots \alpha_{n-j},
\end{equation*}
because $x \in {\vv}_q$. Therefore
\begin{equation}\label{55}
b_{j+1} \ldots b_n > \overline{\alpha_1 \ldots \alpha_{n-j}}.
\end{equation}
Let $(d_i)$ be another expansion of $x$ and let $j$ be the {\it smallest}
positive integer
for which $d_j \not= b_j$. Since $(b_i)$ is greedy, we have $d_j < b_j$
and $j \in \set{1, \ldots , n}$. First we show that $j \in \set{n-k, n}$.
Assume on the contrary that $j \notin \set{n-k,n}$.

First assume that $n-k < j <n$. Then $b_j >0$ and by \eqref{55},
\begin{equation*}
b_{j+1} \ldots b_n 0^{\infty} >
\overline{\alpha_1 \ldots \alpha_{n-j}} \overline{\alpha_{n-j+1} \ldots
\alpha_{k}^-} 0^{\infty}.
\end{equation*}
Since $\alpha_1 \ldots \alpha_{k}^- \alpha_1^{\infty}$ is the smallest
expansion of 1 in base $q$ (see Lemma~\ref{l53}),\newline
$\overline{\alpha_1 \ldots \alpha_k^-}0^{\infty}$ is the greedy expansion of
$\alpha_1 /(q-1) - 1$, and thus
\begin{equation*}
\sum_{i=1}^{\infty} \frac{d_{j+i}}{q^i} =
b_j - d_j + \sum_{i=1}^{\infty} \frac{b_{j+i}}{q^i}  >  \frac{\alpha_1}{q-1},
\end{equation*}
which is impossible.

Next assume that $1 \leq j < n-k$. Rewriting \eqref{55} one gets
\begin{equation*}
\overline{b_{j+1} \ldots b_n} < \alpha_1 \ldots \alpha_{n-j}.
\end{equation*}
If we had
\begin{equation*}
\overline{b_{j+1} \ldots b_{j+k}} = \alpha_1 \ldots \alpha_k,
\end{equation*}
then
\begin{equation*}
\overline{b_{j+k+1} \ldots b_n} < \alpha_{k+1} \ldots \alpha_{n-j}.
\end{equation*}
Hence
\begin{equation*}
b_{j+k+1} b_{j+k+2} \ldots > \overline{\alpha_{k+1} \alpha_{k+2} \ldots} =  \alpha_1 \alpha_2 \ldots.
\end{equation*}
Since in this case $b_{j+k} = \overline{\alpha_k} < \alpha_1$, the last
inequality contradicts the fact that $(b_i)$ is a greedy sequence.
Therefore
\begin{equation*}
\overline{b_{j+1} \ldots b_{j+k}} < \alpha_1 \ldots \alpha_k
\end{equation*}
or equivalently
\begin{equation*}
b_{j+1} \ldots b_{j+k} \geq \overline{\alpha_1 \ldots \alpha_k^-}.
\end{equation*}
Since $n > j+k$ and $b_n >0$, it follows that
\begin{equation*}
b_{j+1} b_{j+2} \ldots > \overline{\alpha_1 \ldots \alpha_k^-}0^{\infty},
\end{equation*}
which leads to the same contradiction as we encountered above.
It remains to investigate what happens if $j \in \set{n-k,n}$.

If $j=n-k$, then it follows from \eqref{55} that
\begin{equation*}
b_{n-k+1} \ldots b_n \geq \overline{\alpha_1 \ldots \alpha_k^-}.
\end{equation*}
Equivalently,
\begin{equation*}
b_{n-k+1} b_{n-k+2} \ldots =
b_{n-k+1} \ldots b_n 0^{\infty} \geq
\overline{\alpha_1 \ldots \alpha_k^-}0^{\infty},
\end{equation*}
and thus
\begin{equation}\label{56}
\sum_{i=1}^{\infty} \frac{d_{n-k+i}}{q^i} \geq
1+ \sum_{i=1}^{\infty} \frac{b_{n-k+i}}{q^i}  \geq \frac{\alpha_1}{q-1},
\end{equation}
where both inequalities in \eqref{56} are equalities if and only if
\begin{equation*}
d_{n-k}=b_{n-k}^-,\,
b_{n-k+1} \ldots b_n = \overline{\alpha_1 \ldots \alpha_k^-}, \,\mbox{ and }
\, d_{n-k+1}d_{n-k+2} \ldots = \alpha_1^{\infty}.
\end{equation*}
Hence $d_{n-k} < b_{n-k}$ is only possible in case $b_{n-k} >0$ and
$b_{n-k+1} \ldots b_n = \overline{\alpha_1 \ldots \alpha_k^-}$.

Finally, if $j=n$, then $d_n=b_n^-$ and $(d_{n+i})$ is one of the
expansions listed in Lemma~\ref{53}.
\end{proof}

\begin{remark}
Fix $q \in {\vv} \setminus {\uuu}$.
By Lemma~\ref{l53} the number 1 has a finite greedy
expansion in base $q$. Hence each element
$q \in {\vv} \setminus {\uuu}$ is algebraic. Since
each
$x \in {\vv}_q \setminus {\uu}_q$ has a finite greedy expansion in
base $q$, it follows that the set ${\vv}_q \setminus {\uu}_q$
consists entirely of algebraic numbers.
\end{remark}

\begin{proof}[Proof of Theorem~\ref{t15}]
Fix $q \in (1,\infty )\setminus {\vv}$. In view of Lemma~\ref{l52} it suffices to
prove that a number $x \in J_q \setminus \set{0}$ with a finite
greedy expansion does {\it not} belong to ${\vv}_q$.

Let $x \in J_q \setminus \set{0}$ be an element with a finite greedy expansion.
Since $q \notin {\vv}$, there exists a positive integer $n$ such that
\begin{equation*}
\overline{\alpha_{n+1} \alpha_{n+2} \ldots} > \alpha_1 \alpha_2 \ldots .
\end{equation*}
Let $m = \max \set{i \in \NN: 1 \leq i \leq n \text{ and } \alpha_i > 0}$. From Proposition~\ref{p23} it follows at once that
\begin{equation*}
\alpha_m > 0 \quad \mbox{and} \quad \overline{\alpha_{m+1} \alpha_{m+2} \ldots} > \alpha_1 \alpha_2 \ldots .
\end{equation*}
Since $(a_i(x))$ ends with $\alpha_1 \alpha_2 \ldots$,
we conclude that $x \notin {\vv}_q.$
\end{proof}

\begin{corollary}\label{c54}
A real number $q>1$ belongs to $\uuu$ if and only if $1$ has a unique {\rm infinite} expansion in base $q$.
\end{corollary}

\begin{proof} Suppose first that $q \in \uuu$. It follows from Theorems~\ref{t11} and ~\ref{t12} that $(\alpha_i(q))$ is the smallest expansion of 1 in base $q$. Since $(\alpha_i(q))$ is by definition its largest infinite expansion in base $q$, it follows that $(\alpha_i(q))$ is the unique infinite expansion of 1 in base $q$. Conversely, if 1 has a unique infinite expansion in base $q$, then $(\alpha_i(q))=(\alpha_i)$ is its smallest expansion because the smallest expansion of 1 is infinite in any base $q>1$. Hence $\overline{\alpha_1 \alpha_2  \ldots}$ is a greedy sequence in base $q$. Applying \eqref{51} (note that $\overline{\alpha_1} = 0$) and Theorem~\ref{t12} we conclude that $q \in \uuu$.
\end{proof}

\section{Proof of Theorems~\ref{t16}, \ref{t17}, \ref{t18} and \ref{t19}}\label{s6}

In this section we will complete our study of the sets $\uu_q$
for numbers $q >1$. The results proved in the preceding sections were mainly
concerned with various properties of these sets for numbers $q \in \vv$.
Now we will use these properties to describe the topological structure
of $\uu_q$ for each number $q >1$.

Since the set $\vv$ is closed, we may write $(1, \infty) \setminus \vv$ as
the union of countably many disjoint open intervals $(q_1,q_2)$: the connected
components of $(1, \infty) \setminus \vv$. In order to determine the endpoints
of these components we recall from \cite{[KL3]} that $\vv \setminus \uuu$ is dense in $\vv$ and
all elements of $\vv \setminus \uuu$ are isolated in
$\vv$. In fact, for each element $q \in \uuu$ there exists a sequence
$(q_m)_{m \geq 1}$ of numbers in $\vv \setminus \uuu$ such that
$q_m \uparrow q$, as can be seen from the proof of Theorem 2.6 in \cite{[KL3]}.

\begin{proposition}\label{p61} \mbox{}
\begin{itemize}
\item[\rm (i)] The set R of right endpoints $q_2$ of the connected components
$(q_1,q_2)$ is given by $R = \vv \setminus \uuu$.
\item[\rm (ii)] The set L of left endpoints $q_1$ of the connected components
$(q_1,q_2)$ is given by $L = {\NN} \cup (\vv \setminus \uu)$.
\end{itemize}
\end{proposition}

\begin{proof}[Proof of Proposition~\ref{p61} (i).]
Note that $\vv \setminus \uuu \subset R$ because the
set $\vv \setminus \uuu$ is discrete. As we have already observed above,
each element $q \in \uuu$ can be approximated arbitrarily closely
by elements of $\vv$ smaller than $q$, and thus
$R = \vv \setminus \uuu$.
\end{proof}

The proof of part~(ii) of Proposition~\ref{p61} requires more
work.  We will first prove a number of technical lemmas.
We recall from Section~\ref{s1} that the notation $q \sim (\alpha_i)$ means that the
quasi-greedy expansion of $1$ in base $q$ is given by
$(\alpha_i)$. For convenience we also write $1 \sim 0^{\infty}$,
and occasionally we refer to $0^{\infty}$ as the quasi-greedy
expansion of the number 1 in base 1.

In Lemmas~\ref{l62} and ~\ref{l63} below, $q_2$ is a fixed (but arbitrary) element of $\vv \setminus \uuu$, and
\begin{equation*}
(\alpha_i)=(\alpha_i(q_2))=
(\alpha_1 \ldots \alpha_k \overline{\alpha_1 \ldots \alpha_k})^{\infty}
\end{equation*}
where $k$ is chosen to be minimal.

\begin{remark}
The minimality of $k$ implies that the least period of $(\alpha_i)$ equals
$2k$. Indeed, if $j$ is the least period of
$(\alpha_i)$, then $\alpha_j = \alpha_{2k} =
\overline{\alpha_k} < \alpha_1$ because $j$ divides $2k$. Hence
$\alpha_1 \ldots \alpha_j^+0^{\infty}$ is an expansion of $1$ in base
$q_2$ which contradicts Lemma~\ref{l53} if $j < 2k$.
\end{remark}

\begin{lemma}\label{l62}
For all $i$ with $0 \leq i < k$, we have
\begin{equation*}
\overline{\alpha_{i+1} \ldots \alpha_k} < \alpha_1 \ldots \alpha_{k-i}.
\end{equation*}
\end{lemma}

\begin{proof} For $i=0$ the inequality follows from the relation
$\overline{\alpha_1}=0 < \alpha_1$. Henceforth assume that
$1 \leq i < k.$ Since $q_2 \in \vv$,
\begin{equation*}
\overline{\alpha_{i+1} \ldots \alpha_k} \leq
\alpha_1 \ldots \alpha_{k-i}.
\end{equation*}
Suppose that
\begin{equation*}
\overline{\alpha_{i+1} \ldots \alpha_k} = \alpha_1 \ldots \alpha_{k-i}.
\end{equation*}
If $k \geq 2i$, then
\begin{equation*}
\alpha_1 \ldots \alpha_{2i}= \alpha_1 \ldots \alpha_i
\overline{\alpha_1 \ldots \alpha_i},
\end{equation*}
and it would follow from Lemma~\ref{l42} that
\begin{equation*}
(\alpha_i)= (\alpha_1 \ldots \alpha_i
\overline{\alpha_1 \ldots \alpha_i})^{\infty},
\end{equation*}
contradicting the minimality of $k$. If $i < k < 2i$, then
\begin{align*}
\overline{\alpha_{i+1} \ldots \alpha_{2i}} & =
\overline{\alpha_{i+1} \ldots \alpha_k} \alpha_1 \ldots \alpha_{2i-k} \\
& = \alpha_1 \ldots \alpha_{k-i} \alpha_1 \ldots \alpha_{2i-k} \\
& \geq \alpha_1 \ldots \alpha_{k-i} \alpha_{k-i+1} \ldots \alpha_i \\
& = \alpha_1 \ldots \alpha_i,
\end{align*}
leading to the same contradiction.
\end{proof}

Let $q_1$ be the largest element of
$\vv \cup \set{1}$ that is smaller than $q_2$.
This element exists because the set $\vv \cup \set{1}$ is closed and
the elements of $\vv \setminus \uuu$ are isolated points of
$\vv \cup \set{1}$.
The next lemma provides the quasi-greedy expansion of 1 in base $q_1$.
\begin{lemma}\label{l63}
$q_1 \sim (\alpha_1 \ldots \alpha_k^-)^{\infty}$.
\end{lemma}

\begin{proof}
Let $q_1 \sim (v_i)$. If $k=1$, then
$q_2 \sim (\alpha_1 0)^{\infty}$, and $(v_i)=(\alpha_1^-)^{\infty}$ because
$q_2$ is the smallest element of $\vv \cap (\alpha_1, \alpha_1 + 1)$.
Hence we may assume that $k \geq 2$. This implies in particular that $q_1 \in \vv$ and
$\lceil q_1 \rceil = \lceil q_2 \rceil$. Observe that
\begin{equation*}
v_1 \ldots v_k  \leq \alpha_1 \ldots \alpha_k.
\end{equation*}
If we had
\begin{equation*}
v_1 \ldots v_k = \alpha_1 \ldots \alpha_k,
\end{equation*}
then
\begin{equation*}
v_{k+1} \ldots v_{2k} \leq
\overline{\alpha_1 \ldots \alpha_k},
\end{equation*}
i.e.,
\begin{equation*}
 \overline{v_{k+1} \ldots v_{2k}} \geq
\alpha_1 \ldots \alpha_k
 = v_1 \ldots v_k
\end{equation*}
and it would follow from Lemma~\ref{l42} that $q_1 = q_2$.
Hence
\begin{equation*}
v_1 \ldots v_k \leq \alpha_1 \ldots \alpha_k^-.
\end{equation*}
It follows from Proposition~\ref{p23} that
$(w_i)= (\alpha_1 \ldots \alpha_k^-)^{\infty}$
is the largest quasi-greedy expansion of 1 in some base $q > 1$ that
starts with
$\alpha_1 \ldots \alpha_k^-.$ Therefore it
suffices to show that the sequence $(w_i)$ satisfies the inequalities \eqref{52}.
Since the sequence $(w_i)$ is periodic with period $k$,
it is sufficient to verify that
\begin{equation}\label{61}
\overline{w_{j+1} w_{j+2} \ldots } \leq w_1 w_2 \ldots
\end{equation}
for all $j$ with $0 \le j < k$. If $j=0$, then \eqref{61} is true because $\overline{w_1}=0 < w_1$;
hence assume that $1 \leq j < k$. Then,
according to the preceding lemma,
\begin{equation*}
\overline{\alpha_{j+1} \ldots \alpha_k} < \alpha_1 \ldots \alpha_{k-j}
\end{equation*}
and
\begin{equation*}
\overline{\alpha_1 \ldots \alpha_j} < \alpha_{k-j+1} \ldots \alpha_k.
\end{equation*}
Hence
\begin{align*}
\overline{w_{j+1} \ldots w_{j+k}} & =
\overline{\alpha_{j+1} \ldots \alpha_k^-}
\overline{\alpha_1 \ldots \alpha_j} \\
& \leq \alpha_1 \ldots \alpha_{k-j} \overline{\alpha_1 \ldots \alpha_j} \\
& < \alpha_1 \ldots \alpha_k,
\end{align*}
so that
\begin{equation*}
\overline{w_{j+1} \ldots w_{j+k}} \leq w_1 \ldots w_k.
\end{equation*}
Since the sequence $(w_{j+i})=w_{j+1} w_{j+2} \ldots$ is also periodic with
period $k$, the inequality \eqref{61} follows.
\end{proof}

We include for completeness the following lemma (see also [KL3]).

\begin{lemma}\label{l64}
Fix $q >1$ and let $(\beta_i)=(b_i(1,q))$ be the greedy expansion of the number $1$
in base $q$. For any positive integer $n$, we have
\begin{equation*}
\beta_{n+1} \beta_{n+2} \ldots \leq \beta_1 \beta_2 \ldots.
\end{equation*}
\end{lemma}

\begin{proof}
Let $n \in \NN$. From \eqref{51} we get that
\begin{equation*}
\beta_{n+1} \beta_{n+2} \ldots < \alpha_1 \alpha_2 \ldots
\leq \beta_1 \beta_2 \ldots,
\end{equation*}
whenever there exists a positive integer $j \leq n$ satisfying
$\beta_j < \beta_1 = \alpha_1$. If such an integer $j$ does not exist,
then either
$(\beta_i) = \alpha_1^{\infty}$ or there exists an integer $j >n$ for which
$\beta_j < \alpha_1$. In both these cases the desired inequality readily
follows as well.
\end{proof}

Now we consider a number $q_1 \in \vv \setminus \uu$. Recall from
Lemmas \ref{l43} and \ref{l53} that the greedy expansion $(\beta_i)$ of 1
in base $q_1$ is finite. Let $\beta_m$ be its last nonzero element.

\begin{lemma}\label{l65}
\mbox{}
\begin{itemize}
\item[\rm (i)]
The least element $q_2$ of $\vv$ that is larger than $q_1$ exists.
Moreover,
\begin{equation*}
q_2 \sim (\beta_1 \ldots \beta_m \overline{\beta_1 \ldots \beta_m})^{\infty}.
\end{equation*}
\item[\rm (ii)]
The greedy expansion of $1$ in base $q_2$ is given by
$(\gamma_i) = \beta_1 \ldots \beta_m \overline{\beta_1 \ldots \beta_m^-}
0^{\infty}.$
\end{itemize}
\end{lemma}

\begin{proof}
(i) First of all, note that
\begin{equation*}
q_1 \sim (\alpha_i)=(\beta_1 \ldots \beta_m^-)^{\infty}.
\end{equation*}
Moreover, $(\beta_1 \ldots \beta_m^-)^{\infty}$ is the largest
quasi-greedy expansion of 1 in some base $q >1$ that starts with
$\beta_1 \ldots \beta_m^-$. Hence, in view of Lemma~\ref{l42},
it suffices to show that the infinite sequence
\begin{equation*}
(w_i)=(\beta_1 \ldots \beta_m \overline{\beta_1 \ldots \beta_m})^{\infty}
\end{equation*}
satisfies the inequalities
\begin{equation}\label{62}
w_{k+1} w_{k+2} \ldots \leq w_1 w_2 \ldots
\end{equation}
and
\begin{equation}\label{63}
\overline{w_{k+1} w_{k+2} \ldots} \leq w_1 w_2 \ldots
\end{equation}
for all $k \geq 0$. Observe that \eqref{62} for $k+m$ is equivalent to
\eqref{63} for $k$ and \eqref{63} for $k+m$ is equivalent to \eqref{62}
for $k$. Since both relations are obvious for
$k=0$, we only need to verify \eqref{62} and \eqref{63} for all $k$ with $1 \leq k < m$.
Fix such an index $k$.

The relation \eqref{63} follows from our assumption that
$q_1 \in \vv$:
\begin{equation*}
\overline{w_{k+1} \ldots w_m} =
\overline{\beta_{k+1} \ldots \beta_m} <
\overline{\alpha_{k+1} \ldots \alpha_m} \leq
\alpha_1 \ldots \alpha_{m-k} =
w_1 \ldots w_{m-k}.
\end{equation*}
Since $1 \leq m-k < m$, we also have
\begin{equation*}
\overline{w_{m-k+1} \ldots w_m} < w_1 \ldots w_k.
\end{equation*}
Using Lemma~\ref{l64} we obtain
\begin{align*}
w_{k+1} \ldots w_{k+m} &= w_{k+1} \ldots w_m \overline{w_1 \ldots w_k} \\
& \leq  w_1 \ldots w_{m-k} \overline{w_1 \ldots w_k} \\
& < w_1 \ldots w_{m-k} w_{m-k+1} \ldots w_m,
\end{align*}
from which \eqref{62} follows.

(ii) The sequence $(\gamma_i)$ is an expansion of 1 in base $q_2$. It remains to show that
\begin{equation}\label{64}
\gamma_{k+1} \gamma_{k+2} \ldots < w_1 w_2 \ldots
\quad \mbox{whenever} \quad \gamma_k < w_1.
\end{equation}
If $1 \leq k < m$, then \eqref{64} follows from
\begin{equation*}
\gamma_{k+1} \ldots \gamma_{k+m} = w_{k+1} \ldots w_{k+m} < w_1 \ldots w_m.
\end{equation*}
If $k=m$, then \eqref{64} follows from $\gamma_{m+1}= \overline{w_1}
=0 < w_1$ (note that $m >1$).\\
If $m < k < 2m$, then
\begin{equation*}
\gamma_{k+1} \ldots \gamma_{2m} = \overline{\beta_{k-m+1} \ldots \beta_m^-}
\leq w_1  \ldots w_{2m-k}.
\end{equation*}
Hence
\begin{equation*}
\gamma_{k+1} \gamma_{k+2} \ldots = \gamma_{k+1} \ldots \gamma_{2m}
0^{\infty} < w_1 w_2 \ldots,
\end{equation*}
because $(w_i)$ is infinite. Finally, if $k \geq 2m$, then $\gamma_{k+1} = 0
< w_1.$
\end{proof}

\begin{proof}[Proof of Proposition~\ref{p61}~(ii).]
It follows from Lemma~\ref{l65} that $\vv \setminus \uu \subset L$.
If $q_2 \sim (n0)^{\infty}$ for some $n \in {\NN}$, then $(n, q_2)$ is
a connected component of $(1, \infty) \setminus \vv$. Hence
${\NN} \subset L$. It remains to show that
$(L \setminus {\NN}) \cap \uu = \varnothing.$

If $(q_1,q_2)$ is a connected component of $(1,\infty) \setminus \vv$ with
$q_2 \sim (\alpha_i)$ and $q_1 \in L \setminus {\NN}$,
then by Proposition~\ref{p61}~(i) and Lemma~\ref{l63}, $q_1 \sim (\alpha_1 \ldots \alpha_k^-)^{\infty}$
for some $k \geq 2$.
Since $\alpha_1 \ldots \alpha_k 0^{\infty}$ is another
expansion of 1 in base $q_1$, we have $q_1 \notin \uu$.
\end{proof}

Recall from Section~\ref{s1} that for $q>1$, $\uu_q'$ and
$\vv_q'$ denote the sets of quasi-greedy expansions in base $q$ of the numbers
$x \in \uu_q$ and $x \in \vv_q$ respectively.

\begin{lemma}\label{l66}  Let $(q_1,q_2)$ be a connected component of
$(1, \infty) \setminus \vv$ and suppose that $q_1 \in \vv \setminus \uu$.
Then
\begin{equation*}
\uu_{q_2}'=\vv_{q_1}'.
\end{equation*}
\end{lemma}

\begin{proof} First of all, note that $\lceil q_1 \rceil = \lceil q_2 \rceil$ because $q_1 \notin \NN$ by assumption. Hence the conjugate bars with respect to $q_1$ and $q_2$ have the same meaning.

Let us write again
\begin{equation*}
q_2 \sim (\alpha_1 \ldots \alpha_k
\overline{\alpha_1 \ldots \alpha_k})^{\infty}
\end{equation*}
where $k$ is chosen to be minimal.
Suppose that a sequence $(c_i) \in \set{0, \ldots , \alpha_1}^{\NN}$
is univoque in base $q_2$, i.e.,
\begin{equation}\label{65}
c_{n+1} c_{n+2} \ldots < (\alpha_1 \ldots \alpha_k
 \overline{\alpha_1 \ldots \alpha_k})^{\infty} \quad \mbox{whenever} \quad
c_n < \alpha_1
\end{equation}
and
\begin{equation}\label{66}
\overline{c_{n+1} c_{n+2} \ldots} < (\alpha_1 \ldots \alpha_k
 \overline{\alpha_1 \ldots \alpha_k})^{\infty} \quad \mbox{whenever} \quad   c_n >0.
\end{equation}
If $c_n < \alpha_1$, then by \eqref{65},
\begin{equation*}
c_{n+1} \ldots c_{n+k} \leq \alpha_1 \ldots \alpha_k.
\end{equation*}
If we had
\begin{equation*}
c_{n+1} \ldots c_{n+k} =  \alpha_1 \ldots \alpha_k,
\end{equation*}
then
\begin{equation*}
c_{n+k+1} c_{n+k+2} \ldots < (\overline{\alpha_1 \ldots \alpha_k}
 \alpha_1 \ldots \alpha_k)^{\infty},
\end{equation*}
and by \eqref{66} (note that in this case $c_{n+k}=\alpha_k >0$),
\begin{equation*}
c_{n+k+1} c_{n+k+2} \ldots >(\overline{\alpha_1 \ldots \alpha_k}
 \alpha_1 \ldots \alpha_k)^{\infty},
\end{equation*}
a contradiction. Hence
\begin{equation*}
c_{n+1} \ldots c_{n+k} \leq \alpha_1 \ldots \alpha_k^-.
\end{equation*}
Note that $c_{n+k} < \alpha_1$ in case of equality.
It follows by induction that
\begin{equation*}
c_{n+1} c_{n+2} \ldots \leq (\alpha_1 \ldots \alpha_k^-)^{\infty}.
\end{equation*}
Since a sequence $(c_i)$ satisfying \eqref{65} and \eqref{66} is
infinite unless $(c_i)=0^{\infty}$, we conclude from
Proposition~\ref{p24} and Lemma~\ref{l63} that  $(c_i)$ is the
quasi-greedy expansion of some $x$ in base $q_1$. Repeating the
above argument for the sequence $\overline{c_1 c_2 \ldots}$, which
is also univoque in base $q_2$, we conclude that $(c_i) \in
\vv_{q_1}'$. The reverse inclusion follows from the fact that the
map $q \mapsto (\alpha_i(q))$ is strictly increasing.
\end{proof}

\begin{lemma}\label{l67}
Let $(q_1,q_2)$ be a connected component of
$(1, \infty) \setminus \vv$ and suppose that $q_1 \in \vv \setminus \uu$.
If $q \in (q_1,q_2]$, then

\begin{itemize}
\item[\rm (i)] $\uu_q' = \vv_{q_1}'$;
\item[\rm (ii)] $\uu_q$ contains isolated points if and only if
$q_1 \in \vv \setminus \uuu$. Moreover, if $q_1 \in \vv \setminus \uuu$, then
each sequence $(a_i) \in \vv_{q_1}' \setminus \uu_{q_1}'$ is the expansion
in base $q$ of an isolated point of $\uu_q$ and each sequence $(c_i) \in \uu_{q_1}'$
is the expansion in base $q$ of an accumulation point of $\uu_q$.
\end{itemize}
\end{lemma}

\begin{proof}

(i) Note that
\begin{equation}\label{67}
\uu_q' \subset \uu_r' \quad \mbox{and} \quad \ \vv_q' \subset \uu_r'
\quad \mbox{if } 1< q < r \mbox{ and } \lceil q \rceil = \lceil r \rceil.
\end{equation}
It follows from Lemma~\ref{l66} and ~\eqref{67} that
$\uu_q'= \vv_{q_1}'$ for all $q \in (q_1, q_2]$.

(ii) We need the following observation (valid for all $q >1$) which is a consequence of
Lemmas \ref{l31} and \ref{l32}:

If $x \in J_q$ has an infinite greedy expansion, then a sequence
$(x_i)$ with elements in $J_q$ converges to $x$ if and only if the greedy
expansion of $x_i$ converges (coordinate-wise) to the greedy expansion of $x$
as $i \to \infty$.
Moreover, $x_i \downarrow 0$ if and only if the greedy
expansion of $x_i$ converges (coordinate-wise)
to the sequence $0^{\infty}$ as $i \to \infty$.

First assume that $q_1 \in \vv \setminus {\uuu}$. Let $x \in
\vv_{q_1} \setminus \uu_{q_1}$ and denote the quasi-greedy
expansion of $x$ in base $q_1$ by $(a_i)$. Since each element in
$\vv_{q_1} \setminus \uu_{q_1}$ is an isolated point of
$\vv_{q_1}$ (see Theorem~\ref{t14}~(ii)), there exists a positive
integer $n$ such that the quasi-greedy expansion in base $q_1$ of
any element in $\vv_{q_1} \setminus \set{x}$ does {\it not} start
with $a_1 \ldots a_n$. Since $(a_i) \in \vv_{q_1}'= \uu_q'$, it
follows from the above observation that the sequence $(a_i)$ is
the unique expansion in base $q$ of an isolated point of $\uu_q$.
If $x \in \uu_{q_1}$, then there exists a sequence of numbers
$(x_i)$ with $x_i \in \vv_{q_1} \setminus \uu_{q_1}$ such that the
quasi-greedy expansions of the numbers $x_i$ converge to the
unique expansion of $x$, as can be seen from the proof of
Theorem~\ref{t14}~(iib) (which in turn relies on the proof of
Theorem~\ref{t13}~(iib)). Hence the unique expansion of $x$ in base
$q_1$ is the unique expansion in base $q$ of an accumulation point
of $\uu_q$.

Next assume that $q_1 \in {\uuu} \setminus \uu$. It follows from Theorem~\ref{t13}
that the set $\uu_{q_1}$ has no isolated points.
Hence for each $x \in \uu_{q_1}$ there exists a sequence of
numbers $(x_i)$ with $x_i \in \uu_{q_1} \setminus \set{x}$  such
that $x_i \to x$. In view of the above observation, the unique
expansions of the numbers $x_i$ converge to the unique expansion
of $x$. Therefore the unique expansion of $x$ in base $q_1$ is the
unique expansion in base $q$ of an accumulation point of $\uu_q$.
If $x \in \vv_{q_1} \setminus \uu_{q_1}= \overline{\uu_{q_1}}
\setminus \uu_{q_1}$, then there exists a sequence $(x_i)$ of
numbers in $\uu_{q_1}$ such that the unique expansions of the
numbers $x_i$ converge to the quasi-greedy expansion $(a_i)$ of
$x$, as follows from the proof of Theorem~\ref{t13}~(i).
Hence, also in this case, $(a_i)$ is the unique expansion in base
$q$ of an accumulation point of $\uu_q$. Since $\uu_q' =
\vv_{q_1}'$, this completes the proof.
\end{proof}

\begin{lemma}\label{l68}
Let $(q_1,q_2)$ be a connected component of
$(1, \infty) \setminus \vv$ and suppose that $q_1 \in {\NN}$.
If $q \in (q_1, q_2]$, then $\uu_q' = \uu_{q_2}'$
and $\uu_q$ contains isolated points if and only if
$q_1 \in \set{1,2}$.
\end{lemma}

\begin{proof} Note that if $q_1 = n \in {\NN}$, then
$q_2 \sim (n 0) ^{\infty}$.
Suppose that $q \in (n, q_2]$.
We leave the verification of the following
statements to the reader.\\
A sequence $(a_i) \in \set{0, \ldots , n}^{\NN}$
belongs to $\uu_q'$ if and only if for all $j \in \NN$,
\begin{equation*}
a_j < n \Longrightarrow a_{j+1} < n
\end{equation*}
and
\begin{equation*}
a_j > 0 \Longrightarrow a_{j+1} > 0.
\end{equation*}
In particular we see that
\begin{equation*}
\uu_q' = \uu_{q_2}'.
\end{equation*}
If $n=1$, then ${\uu}_q' = \set{0^{\infty},1^{\infty}}$. If $n=2$, then
\begin{equation*}
\uu_q' = \set{0^{\infty}, 2^{\infty}} \cup \bigcup_{n=0}^{\infty}
\set{0^n1^{\infty}, 2^n 1^{\infty}}.
\end{equation*}
Hence, if $n=2$, then $\uu_q$ is countable and all elements of $\uu_q$ are
isolated, except for its endpoints.
If $n \geq 3$, then $\uu_q$ has no isolated points.
\end{proof}

\begin{proposition}\label{p69}
Let $q >1$ be a real number.
\begin{itemize}
\item[\rm (i)] If $q \in \uuu$, then $q$ is neither stable from below
nor stable from above.
\item[\rm (ii)] If $q \in \vv \setminus \uuu$, then $q$ is stable from below,
but not stable from above.
\item[\rm (iii)] If $q \notin \vv$, then $q$ is stable.
\end{itemize}
\end{proposition}

\begin{proof}
(i) As mentioned at the beginning of this section, if $q \in \uuu$,
then there exists a sequence $(q_m)_{m \geq 1}$ with
numbers $q_m \in \vv \setminus \uuu$, such that $q_m \uparrow q$.
Since
\begin{equation*}
\uu_{q_m}' \subsetneq \vv_{q_m}' \subset \uu_q',
\end{equation*}
$q$ is not stable from below. If $q \in \uuu \setminus {\NN}$, then
$q$ is not stable from above because
\begin{equation}\label{68}
\uu_q' \subsetneq \vv_q' \subset \uu_s'
\end{equation}
for any $s \in (q, \lceil q \rceil]$. If $q \in \set{2,3,\ldots}$, then $q$ is
not stable from above because the sequence $q^{\infty}$ belongs to $\uu_s'
\setminus \uu_q'$ for any
$s > q$.

(ii) and (iii) If $q \notin \uuu$, then $q \in (q_1, q_2]$, where
$(q_1, q_2)$ is a connected component of $(1, \infty) \setminus \vv$.
From Proposition~\ref{p61} and Lemmas~\ref{l67} and \ref{l68} we conclude that $q$ is stable
from below. Note that $q=q_2$ if and only if $q \in \vv \setminus \uuu$.
Hence, if $q \notin \vv$, then $q$ is also stable from above.
If $q \in \vv \setminus \uuu$, then $q$ is not stable from above because \eqref{68} holds
for any $s \in (q, \lceil q \rceil]$.
\end{proof}

\begin{remark} The issue of (non-)stability of bases is further explored in \cite{[DV2]}.
\end{remark}

\begin{proof}[Proof of Theorem~\ref{t16}.]
(i) This is the content of Proposition~\ref{p61}.

(ii) If $q \in \set{2,3, \ldots }$, then $\uu_q \subsetneq \overline{\uu_q}
=[0,1].$ Hence, neither $\uu_q$ nor $\overline{\uu_q}$ is a Cantor set.

(iii) and (iv) If $q \notin {\NN}$, then $\uu_q$ is nowhere dense
by a remark following the statement of Theorem~\ref{t15} in
Section~\ref{s1}. Hence, if $q \notin {\NN}$, then $\uu_q$ is a
Cantor set if and only if $\uu_q$ is closed and does not contain
isolated points.

If $q \in \uuu \setminus {\NN}$, then by Theorem~\ref{t13}  the
set $\uu_q$ is not closed, and $\overline{\uu_q}$ has no isolated
points, from which part~(iii) follows.

Finally, let $q \in (q_1, q_2]$, where $(q_1,q_2)$ is a connected component
of $(1, \infty) \setminus \vv$. Since $q \notin \uuu$, the set $\uu_q$ is
closed. It follows from Lemmas \ref{l67} and
\ref{l68} that $\uu_q$ is a Cantor set if and only if
$q_1 \in \set{3,4, \ldots } \cup (\uuu \setminus \uu)$. A quick examination
of the proof of these lemmas yields the last statement of part~(iv).
\end{proof}

\begin{proof}[Proof of Theorem~\ref{t17}.]
The statements of this theorem readily follow from Proposition~\ref{p61}, Lemmas \ref{l67} and \ref{l68}, and
Proposition~\ref{p69}.
\end{proof}

\begin{proof}[Proof of Theorem~\ref{t18}.]
  Fix $q \in (1, \infty) \setminus \uuu$. Then
$q \in (q_1, q_2]$, where $(q_1, q_2)$ is a connected component of
$(1, \infty) \setminus \vv$. Let us write
\begin{equation*}
q_2 \sim (\alpha_i) =
(\alpha_1 \ldots \alpha_k \overline{\alpha_1 \ldots \alpha_k})^{\infty},
\end{equation*}
where $k$ is minimal. Let
\begin{equation*}
\mathcal{F} = \set{ja_1 \ldots a_k \in  \set{0, \ldots , \alpha_1}^{k+1} \,: \, j < \alpha_1 \quad \mbox{and}
\quad a_1 \ldots a_k \geq \alpha_1 \ldots \alpha_k}.
\end{equation*}

It follows from Lemma~\ref{l67} (i) and
the proof of Lemmas \ref{l66} and \ref{l68} that a
sequence $(c_i) \in \set{0, \ldots, \alpha_1}^{\NN}$ belongs to
$\uu_q'$ if and only if
$c_j \ldots c_{j+k} \notin \mathcal{F}$ and
$\overline{c_j \ldots c_{j+k}} \notin \mathcal{F}$ for all $j \geq 1$.
Therefore, $\uu_q'$ is a subshift of finite type. Note that any subshift
$S \subset \set{0, \ldots , \alpha_1}^{\NN}$ is closed in the topology of coordinate-wise convergence.
It remains to show that $\uu_q'$ is not closed if $q \in \uuu$.
 It follows from the proof of Theorem~\ref{t13} (i)
that for each $q \in \uuu$ and $x \in \vv_q \setminus
\uu_{q}$, there exists a sequence $(x_i)$ of numbers in $\uu_q$
such that the unique expansions of the numbers $x_i$ converge to the
quasi-greedy expansion of $x$. Hence the set $\uu_q'$ is not closed if $q \in \uuu$.
\end{proof}

\begin{proof}[Proof of Theorem~\ref{t19}.]
(i) Note that $G \sim (10)^{\infty}$. It follows from
the proof of Lemma~\ref{l68} that $\uu_q' = \set{0^{\infty}, 1^{\infty}}$
for all $q \in (1, G]$.

(ii) Due to the properties of the set $\vv \setminus \uuu$ which were mentioned
at the beginning of this section, we may write
\begin{equation*}
\vv \cap (1,q') = \set{q_n : n \in {\NN}} \quad \mbox{and} \quad
\vv \cap (2, q'')= \set{r_n : n \in {\NN}},
\end{equation*}
 where the $q_n$'s and the $r_n$'s are written in strictly increasing order. Note that
$q_1 \sim (10)^{\infty}$ and $r_1 \sim (20)^{\infty}$. Moreover,
$q_n \uparrow q'$ and $r_n \uparrow q''$. Thanks to \eqref{67} we only need to verify that the sets $\uu_{q_n}$ and
$\uu_{r_n}$ are countable for each $n \in {\NN}$.
It follows from the proof of Lemma~\ref{l68} that $\uu_{q_1}$ and $\uu_{r_1}$
are countable. Now suppose that $\uu_{q_n}$ is countable for some $n \geq 1$.
By Lemma~\ref{l66} we have
\begin{equation*}
\uu_{q_{n+1}}' = \vv_{q_n}' = \uu_{q_n}'
\cup (\vv_{q_n}' \setminus \uu_{q_n}').
\end{equation*}
According to Theorem~\ref{t14} (ii) the set $\vv_{q_n}' \setminus
\uu_{q_n}'$ is countable, whence $\uu_{q_{n+1}}$ is countable as
well. It follows by induction that $\uu_{q_n}$ is countable for
each $n \in {\NN}$. Similarly, $\uu_{r_n}$ is countable for each
$n \in {\NN}$. 

(iii) It follows from Theorem~\ref{t13} that $\vert \uu_{q'} \vert
=2^{\aleph_0}$ and $\vert \uu_{q''} \vert =2^{\aleph_0}$. The
relation \eqref{67} yields that $\vert \uu_q \vert = 2^{\aleph_0}$
for all $q \in [q', 2] \cup [q'',3]$. If $q
>3$, then $\vert \uu_q \vert = 2^{\aleph_0}$ because $\uu_q'$ contains all sequences consisting of merely
ones and twos.
\end{proof}

We conclude this paper with an example and some remarks.

\begin{example}
For any given $k \in {\NN}$ define the numbers $p(k)$ and $q(k)$ by setting
\begin{equation*}
p(k) \sim (1^{k-1}0)^{\infty} \quad \mbox{and} \quad
q(k) \sim (1^k0^k)^{\infty}.
\end{equation*}
It follows from Lemma~\ref{l65} and Theorem~\ref{t17} that the
sets $(p(k), q(k)]$ are maximal stability intervals. Moreover,
it follows from the proof of Theorem~\ref{t18} that a sequence
$(c_i) \in \set{0,1}^{\NN}$ belongs to $\uu_q'$ for $q \in
(p(k), q(k)]$ if and only if a zero is never followed by $k$
consecutive ones and a one is never followed by $k$ consecutive
zeros. This result was first established by Dar\'oczy and K\'atai
in \cite{[DK1]}, using a different approach.

From Lemma~\ref{l65} we know that the least element of $\vv$ larger than $q(k)$ is given by $r(k)$,
where
\begin{equation*}
r(k) \sim (1^k 0^{k-1} 1 0^k 1^{k-1} 0 )^{\infty}.
\end{equation*}
Therefore the sets
$(q(k), r(k)]$ are also maximal stability intervals. If
$q \in (q(k), r(k)]$, then $\uu_q$ is not a Cantor set because
$q(k) \in \vv \setminus \uuu$. Moreover, a number $x \in \uu_q$ is an isolated point of $\uu_q$ if and only if its unique expansion belongs to
$\vv_{q(k)}' \setminus \uu_{q(k)}'$.

Finally, let $k\ge3$, and let $s(k)$ be the least element of $\vv$ larger than $r(k)$. Fix $q \in (r(k), s(k)]$ and let $(c_i)$ be the unique expansion of a number $x \in \uu_q$. Combining Lemma~\ref{l67} with the final remark below, one can ``decompose'' $\uu_q$ as follows:
\begin{itemize}
\item $x$ is an isolated point of $\uu_q$ if and only if  $(c_i) \in \vv_{r(k)}' \setminus \uu_{r(k)}'$.
\item $x$ is an accumulation point of $\uu_q$ if and only if $(c_i) \in \uu_{r(k)}'$.
\item $x$ is a condensation point of $\uu_q$ if and only if $(c_i) \in \uu_{q(k)}'$.
\end{itemize}
\end{example}

\begin{remarks}
\mbox{}
\begin{itemize}
\item Let us now consider the set $L'$ of left endpoints and the set $R'$ of
right endpoints of the connected components of
$(1, \infty) \setminus \uuu$. We will show that
\begin{equation*}
L' = {\NN} \cup (\uuu \setminus \uu) \quad \mbox{and}
\quad  R' \subsetneq \uu \setminus{\NN}.
\end{equation*}
Fix a number
$q \in (1, \infty) \setminus \uuu$. Let $q_1$ be the least element of
$\vv$ satisfying $q_1\ge q$. Since $q_1 \in \vv \setminus \uuu$,
it is both a left endpoint and a right endpoint
of a connected component of
$(1, \infty) \setminus \vv$. Hence there exists a sequence
$q_1 < q_2 < \cdots$ of numbers in $\vv \setminus \uuu$ satisfying
\begin{equation*}
 (q_i, q_{i+1}) \cap \vv = \varnothing \quad \mbox{for all} \quad i\geq 1.
\end{equation*}
Let $\beta_m$ be the last nonzero element of the greedy expansion $(\beta_i)$
of the number 1 in base $q_1$. We define a sequence $(c_i)$ by induction as
follows. First, set
\begin{equation*}
c_1 \ldots c_m = \beta_1 \ldots \beta_m.
\end{equation*}
Then, if $c_1 \ldots c_{2^Nm}$ is already defined for some nonnegative integer
$N$, set
\begin{equation*}
c_{2^Nm +1} \ldots c_{2^{N+1}m -1} = \overline{c_1 \ldots c_{2^Nm-1}}
\quad \mbox{and} \quad c_{2^{N+1}m} = \overline{c_{2^Nm}} +1.
\end{equation*}
It follows from Lemma~\ref{l65} that the greedy expansion of 1 in base $q_n$
is given by $c_1 \ldots c_{2^{n-1}m} 0^{\infty}$ $(n \in \NN)$. Hence $(c_i)$ is an
expansion of 1 in base $q^*$ where
\begin{equation*}
q^* = \lim_{n \to \infty} q_n.
\end{equation*}
Note that $q^* \in \vv$ because $\vv$ is closed. The number $q^*$ cannot belong to $\vv \setminus \uuu$ because this set is
discrete, and it cannot belong to $\uuu \setminus \uu$ because $\uu$ is closed from above (\cite{[KL3]}). Hence $q^* \in \uu$ and
$R' \subset \uu \setminus \NN$. The set $R'$ is a proper subset of $\uu \setminus \NN$ because the latter set is uncountable.

Now let $\ell_1$ be the largest element of $\vv \cup \set{1}$
that is smaller than $q_1$.
Let us also write $\ell_1 \sim (\eta_i)$ and $q_1 \sim (\alpha_i)$.
It follows from
Lemma~\ref{l63} and the remark preceding Lemma~\ref{l62} that $(\eta_i)$
has a smaller period than the least period of $(\alpha_i)$. Hence there exists a finite set of
numbers $\ell_k < \cdots < \ell_1$ in $\vv \cup \set{1}$,
such that for $i$ with $1 \leq i < k$,
\begin{equation*}
(\ell_{i+1}, \ell_i) \cap \vv = \varnothing,
\end{equation*}
and such that $\ell_k$ is a left endpoint of a connected component of
$(1, \infty) \setminus \vv$, but not a right endpoint. This means that
\begin{equation*}
\ell_k \in {\NN} \cup (\uuu \setminus \uu) \quad \mbox{and} \quad
(\ell_k, q) \cap \uuu = \varnothing.
\end{equation*}
Hence $\ell_k \in {\uuu}\cup\set{1}$ and therefore $\ell_k \in L'$.
We may thus conclude that
$L' \subset {\NN} \cup (\uuu \setminus \uu)$. On the other hand,
$L \cap \uuu \subset L'$ because $\uuu \subset \vv$. Taking into account that $1 \in L'$, we deduce from Proposition~\ref{p61}
that $L' = {\NN}\cup (\uuu \setminus \uu)$.
\item The analysis of the preceding remark enables us also to determine for
each $n \in {\NN}$ the least element $q^{(n)}$ of the set
$\uu \cap (n, n+1)$:

Fix $n \in {\NN}$ and let $q$ be the least element of
$\vv \cap (n, n+1)$. Then $q \sim (n0)^{\infty}$ and the greedy expansion
$(\beta_i)$ of 1 in base $q$ is given by $n10^{\infty}$. The sequence $(c_i)$
constructed in the preceding remark with $m=2$ and $c_1c_2 = n1$ is the
unique expansion of the number 1 in base $q^{(n)}$.
\item In \cite{[KL2]} it was shown that for each $n \in {\NN}$, there exists
a smallest number $r^{(n)} >1$ for which exactly one sequence $(c_i) \in \set{0, 1, \ldots, n}^{\NN}$
satisfies the equality
\begin{equation*}
\sum_{i=1}^{\infty} c_i (r^{(n)})^{-i}=1.
\end{equation*}
Although this might appear as an equivalent definition of the numbers
$q^{(n)}$, there is a subtle difference: it is not required that $r^{(n)}$ belongs to $(n,n+1)$. It can be seen from
the results in [KL2] that $r^{(n)} = q^{(n)}$ if and only if
$n \in \set{1,2}$. If $n > 2$, then $r^{(n)} < n < q^{(n)}$.
\item Finally, we determine the condensation points of $\uu_q$ for $q>1$.
If $q \in \uuu$, then each element of $\uu_q$ is a
condensation point of $\uu_q$ because the set $\overline{\uu_q}$ is perfect and $\overline{\uu_q} \setminus \uu_q$ is countable.
Henceforth assume that $q \in (1, \infty) \setminus \uuu$. Let $(p_1,p_2)$ be the connected component of $(1, \infty) \setminus \uuu$
containing $q$, and let $r$ be the least element of $(p_1,p_2) \cap \vv$. Since $\uu_q$ is a closed set, it can be written uniquely as a disjoint union of a countable set $C$ and a perfect set $P$ consisting precisely of the condensation points of $\uu_q$. We claim that
\begin{equation*}
P=\set{\sum_{i=1}^{\infty} c_i q^{-i} : (c_i) \in \uu_r'},
\end{equation*}
except when $p_1 \in \set{1,2}$, in which case $P=\varnothing$. Indeed, if $p_1 \notin \set{1,2}$, then $\uu_r$ is a Cantor set by Theorem~\ref{t16} (iv). It follows in particular that the set $\uu_r$ is perfect and hence consists entirely of condensation points. Note that $\uu_r' \subset \uu_q'$ because $(p_1,r]$ is a stability interval and $\lceil r \rceil = \lceil q \rceil$. If $q \in (p_1,r]$, then $\uu_r' = \uu_q'$. If $q \in (r,p_2)$, then $[r,q) \cap \vv$ is a finite subset of $\vv \setminus \uuu$, as follows from the first remark above. Moreover, if we write $[r,q) \cap \vv = \set{r_1, \ldots,r_m}$ where $r_1 < \cdots < r_m$, then by applying Lemma~\ref{l67}~(i) ($m$ times), we get that
\begin{equation*}
\uu_q' = \uu_r' \cup \bigcup_{\ell=1}^m (\vv_{r_{\ell}}' \setminus \uu_{r_{\ell}}').
\end{equation*}
Hence $\uu_q' \setminus \uu_r'$ is countable by Theorem~\ref{t14}~(ii).

Now let $x \in \uu_q$ and let $(c_i)$ be its unique expansion in base $q$. Suppose first that $(c_i)$ belongs to $\uu_r'$. Let $W$ be an arbitrary neighborhood of $x$. There exists an index $N$ such that each univoque sequence in base $q$ starting with the block $c_1 \ldots c_N$ is the unique expansion in base $q$ of a number in $W$. Applying Lemmas~\ref{l31} and ~\ref{l32} to greedy expansions in base $r$ and using the fact that $\uu_r$ is perfect and $\uu_r' \subset \uu_q'$, we conclude that $x$ is a condensation point of $\uu_q$. If $(c_i)$ does not belong to $\uu_r'$, then there exists an index $N$ such that no sequence in $\uu_r'$ starts with $c_1 \ldots c_N$ because $\uu_r$ is closed. Applying Lemmas~\ref{l31} and ~\ref{l32} to greedy expansions in base $q$ and using the fact that $\uu_q' \setminus \uu_r'$ is countable, we conclude that $x$ is not a condensation point of $\uu_q$. If $p_1 \in \set{1,2}$, then $\uu_q$ is countable by Theorem~\ref{t19} and thus $P=\varnothing$.
\end{itemize}
\end{remarks}

\noindent
{\it Acknowledgement.} We thank the referee for his/her helpful remarks. The first author is grateful to the team of the
Erwin Schr\"odinger Institute for Mathematical Physics at Vienna, and in
particular to Klaus Schmidt, for their hospitality and financial support
during the spring of 2006.

\end{document}